\documentclass [a4paper,12pt]{article}

\usepackage{amsthm}
\usepackage{amsmath,amssymb,latexsym,amsfonts,mathrsfs}
\usepackage[dvipdfmx]{graphicx}

\usepackage{color}

\usepackage{amsrefs}
\usepackage{mathtools}
\usepackage{fancyhdr}
\usepackage[top=30truemm, bottom=30truemm, left=25truemm, right=25truemm]{geometry}
\usepackage{physics}

\allowdisplaybreaks[4]

\newcommand{\ep}{\varepsilon}

\newcommand{\SCR}[1]{{\mathscr #1}}
\newcommand{\CAL}[1]{{\cal #1}}
\newcommand{\J}[1]{\left\langle #1 \right\rangle}
\newcommand{\D}[1]{{\mathscr D}( #1 )}

\theoremstyle{definition}
\newtheorem{Thm}{{\bf Theorem}}[section]
\newtheorem{Lem}[Thm]{{\bf Lemma}}
\newtheorem{Prop}[Thm]{{\bf Proposition}}

\newtheorem{Assu}[Thm]{{\bf Assumption}}

\newtheorem{Rem}[Thm]{{\bf Remark}}

\usepackage{bm}
\usepackage[dvipdfmx]{hyperref}
\hypersetup{
    colorlinks=true, 
    linktoc= page,   
}

\def\({\left(}
\def\){\right)}
\def\<{\left\langle}
\def\>{\right\rangle}
\def\le{\leqslant}

\newcommand{\C}{{\bf C}}
\newcommand{\R}{{\bf R}}
\newcommand{\Z}{\mathbb{Z}}

\newcommand{\F}{\mathcal{F}}

\renewcommand{\Im}{\operatorname{Im}}
\renewcommand{\Re}{\operatorname{Re}}
\newcommand{\I}{\infty}
\newcommand{\wha}[1]{\widehat{#1}}

\newcommand{\Lebn}[2]{\left\lVert #1 \right\rVert_{#2}}

\newcommand{\hSobn}[2]{\left\lVert #1 \right\rVert_{\dot{H}^{#2}}}

\begin{document}

\begin{flushleft}
{\Large \bf Long-range scattering for a critical homogeneous type nonlinear Schr\"odinger equation with time-decaying harmonic potentials}
\end{flushleft}

\begin{flushleft}
{\large Masaki KAWAMOTO}\\
{Department of Engineering for Production, Graduate School of Science and Engineering, Ehime University, 3 Bunkyo-cho Matsuyama, Ehime, 790-8577, Japan }\\
Email: {kawamoto.masaki.zs@ehime-u.ac.jp}
\end{flushleft}
\begin{flushleft}
{\large Hayato MIYAZAKI}\\
{Teacher Training Courses, Faculty of Education, Kagawa University, Takamatsu, Kagawa 760-8522, Japan} \\
{E-mail: miyazaki.hayato@kagawa-u.ac.jp}
\end{flushleft}

\begin{abstract}
This paper is concerned with the final state problem for the homogeneous type nonlinear Schr\"odinger equation with time-decaying harmonic potentials. 
The nonlinearity has the critical order and is not necessarily the form of a polynomial.
In the case of the gauge-invariant power-type nonlinearity, the first author proves that the equation admits a nontrivial solution that behaves like a free solution with a logarithmic phase correction in \cite{Ka2}.
In this paper, we extend his result into the case with the general homogeneous nonlinearity by the technique due to the Fourier series expansion introduced by Masaki and the second author \cite{MM2}. 
To adapt the argument in the aforementioned paper, we develop a factorization identity for the propagator and require a little stronger decay condition for the Fourier coefficients arising from the harmonic potential.
Moreover, in two or three dimensions, we improve the regularity condition of the final data in \cite{MM2, MMU}. 
\end{abstract}

\begin{flushleft}

\end{flushleft}

\section{Introduction}
In this paper, we consider the nonlinear Schr\"{o}dinger equation with time-decaying harmonic potentials
\begin{align}\label{eq1} \tag{NLS}
	i \partial_t u - H_0(t) u =  F(u),
\end{align}
where $(t,x) \in {\bf R}^{+} \times {\bf R}^d$, $d \leq 3$, and $u = u(t,x)$ is a complex-valued unknown function.
The operator $H_0(t)$ is defined by
\begin{align*}
	H_0 (t) = - \frac{1}{2}\Delta + \sigma (t) \frac{|x|^2}{2}, \quad \sigma (t) \in \R.
\end{align*}
The coefficient of the potential $\sigma (t)$ satisfies a time-decay condition (see Assumption \ref{A1}). 
Suppose that the nonlinearity $F$ is homogeneous of degree $1+ 2/d(1-\lambda)$; that is, it satisfies 
\begin{align}
	F(\alpha z) = \alpha^{1+p_c} F(z), \quad \lambda \in [0,1/2), \quad p_c \coloneqq \frac{2}{d(1-\lambda)}  \label{noho}
\end{align}
for any $\alpha>0$ and $z \in \C$.
The aim of the paper is to determine the asymptotic behavior of nontrivial solutions to \eqref{eq1} with a general homogeneous nonlinearity.
More specifically, we give a sufficient condition on the shape of the nonlinearity $F$ for \eqref{eq1} to have a short- or long-range scattering solution.
Hence we focus on dealing with the final state problem which is usually discussed simper than the corresponding initial value problem (e.g., \cite{HN98, K21, KM21}).

As for \eqref{eq1} without the potential
\begin{align}\label{eq:f}
	i \partial_t u + \frac12 \Delta u =  F(u),
\end{align}
the asymptotic behavior of solutions for large time has been much studied by many authors. 
We remark that our equation \eqref{eq1} includes the free equation \eqref{eq:f} when $\sigma(t) \equiv 0$ (i.e., $\lambda = 0$).   
The typical nonlinearity is the gauge-invariant power-type nonlinearity 
\begin{align}
	F(u) = \mu |u|^{p} u, \quad \mu \in \R \setminus \{0\}, \quad p>0. \label{gauge:1} 
\end{align}
In the case of \eqref{eq:f} with \eqref{gauge:1}, it is known that the exponent $p_c = 2/d$ is a threshold in view of the long-time behavior of solutions.
It is because the equation admits a nontrivial solution which behaves like the asymptotics of the free solution 
\begin{align}
	(it)^{-\frac{d}{2}}  e^{i \frac{|x|^2}{2t}} \wha{u_{+}}\(\frac{x}{t}\) \label{eq:af}
\end{align}
for large time, when $p> p_c$. 
However, if $p \leq p_c$, then \eqref{eq:f} has no solutions which behave like \eqref{eq:af} in $L^2$ (e.g., \cite{S74, B84, TY84}).
In particular, when $p= p_c$, there is the solution that behaves like a free solution with a logarithmic phase correction
\begin{align}
	(it)^{-\frac{d}{2}}  e^{i \frac{|x|^2}{2t}} \wha{u_{+}}\(\frac{x}{t}\) \exp \( - \mu i \left| \wha{u_{+}}\(\frac{x}{t}\) \right|^{\frac{2}{d}} \log t \) \label{eq:ma}
\end{align}
for a suitable given \emph{final data} $u_+$ as $t \to \infty$, where $\wha{u_+}$ denotes the Fourier transform of $u_+$ (cf. \cite{O91, GO93, HN98}).
Here we say that the nonlinearity is \emph{short-range} if \eqref{eq:f} admits a nontrivial solution which asymptotically behaves like \eqref{eq:af} for large time. 
Also the nonlinearity is said to be \emph{long-range} if \eqref{eq:f} admits a nontrivial solution which asymptotically behaves like \eqref{eq:ma} for large time with a suitable $\mu \in \R \setminus \{0\}$.
We further observe that in $p = p_c$, 
the asymptotic behavior of the solutions depends on the shape of the nonlinearity. 
For instance, when $d=2$, $F(u) = \mu |u| u + \lambda_1 u^2 + \lambda_2 \overline{u}^2$ with $\lambda_j \in \C$ is short-range if $\mu = 0$ and long-range if $\mu \neq 0$ (cf. \cite{HNST04, HNW, ST04}). 
Eventually, via the Fourier series expansion, Masaki and the second author \cite{MM2} treats general nonlinearity satisfying \eqref{noho}, namely, including the non-gauge-invariant nonlinearity, and  
prove that if $g_0=0$ and $g_1 \in \R$, then \eqref{eq:f} admits a nontrivial solution 
which asymptotically behaves like \eqref{eq:ma} with $\mu = g_1$ by introducing a decomposition of the nonlinearity
\begin{align}
	F(u) = g_0 |u|^{1+p_c} + g_1 |u|^{p_c}u + \sum_{n\neq 0,1} g_n |u|^{1+p_c-n}u^n \label{eq:exp}
\end{align}
with the coefficients
\begin{align}
	g_n \coloneqq \frac1{2\pi} \int_0^{2\pi} F(e^{i\theta}) e^{-in\theta} d\theta \label{eq:gn}
\end{align}
under a summability assumption on $\{g_n\}_n$ (see also \cite{MMU}). 
In particular, if $g_0=g_1=0$ then there exists an asymptotically free solution.
We note that if $g_0 \neq 0$, then there are no solutions which behaves like a free solution (cf. \cite{S05, ST06, MM19, MS20}).
Also, we only consider the case where the given final data $u_+$ has a very small low-frequency part, because if the data has a non-negligible low-frequency part, then there appears other kind of asymptotic behavior (see \cite{HN02, HN04-2, HN08, HN11, HN15, N15, NS11}).

Regarding \eqref{eq1}, under a time decay condition for $\sigma(t)$, Kawamoto \cite{Ka4}, and Ishida and Kawamoto \cite{IS1, IS2} discover that the quantum particle governed by the energy $H_0(t)$ is decelerated by the harmonic oscillators with velocity $v = \CAL{O} (|t|^{- \lambda})$ as $t \to \infty $, but which is not trapped in studies on the corresponding linear problem. 
This phenomenon leads us to expect that we can consider the scattering problem for \eqref{eq1} and a kind of new observation different from \eqref{eq:f} appears.
Later on, Kawamoto-Muramatsu \cite{KM21} studies scattering solutions on the initial value problem for \eqref{eq1} with mostly power-like nonlinearities. 
In \cite{KM21}, they find that the threshold $p_c$ changes from $2/d$ to $2/d(1- \lambda)$ and specifies the asymptotics of the solutions for large time, corresponding to \eqref{eq:ma}.
Also, when $p=p_c$, Kawamoto \cite{KM21} proves that \eqref{eq1} with \eqref{gauge:1} admits a nontrivial solution which behaves like the asymptotics for large time.
In this paper, we aim to extend the result of \cite{KM21} into general nonlinearities including non-gauge-invariant types, namely, satisfying \eqref{noho}.

In what follows, in order to investigate large time behavior of solutions to \eqref{eq1}, we assume the following time decay condition on the coefficient $\sigma (t)$:
\begin{Assu} \label{A1}
Let $\zeta _1 (t)$ and $\zeta _2 (t)$ be the fundamental solutions to 
\begin{align*}
\zeta _j ''(t) + \sigma (t) \zeta _j (t) =0, \quad
\begin{cases}
\zeta _1 (0) = 1, \\
\zeta _1 ' (0) =0,
\end{cases}
\quad
\begin{cases}
\zeta _2 (0) = 0, \\
\zeta _2 ' (0) =1.
\end{cases}
\end{align*} 
Then, there exists $c>0$, $r_0 >0$, $\sigma_0$, $c_{1}^{(k)}$, $c_{2}^{(k)} \notin \{0, \infty \}$ and $c_{3} \in {\bf R}$, and $\lambda \in [0, 1/2 )$ such that  $\sigma \in C^1([r_0, \I), \R)$ with $\lim_{t \to \I} t^{3} \sigma' (t) = \sigma_0$, 
$|\zeta _2 (t)| \geq c$ for all $t > r_0$, and for all $k \in \{0, 1, 2\}$, 
\begin{align}\label{K13}
	\lim_{t \to \infty} \frac{|\zeta_1^{(k)} (t)|}{t^{\lambda-k} } = c_{1}^{(k)}, \quad 
	\lim_{t \to \infty} \frac{|\zeta_2^{(k)} (t)|}{t^{1-\lambda-k}} = c_{2}^{(k)}, \quad 
	\lim_{t \to \infty} \frac{\left| \zeta _2(t) -c_{2}^{(0)} t^{1- \lambda} \right|}{t^{\lambda}} = c_{3}
\end{align}
hold. 
\end{Assu}
In the followings, we further assume $\zeta_1(t),\, \zeta_2(t) \geq 0$ for any $t >r_0$ for simplicity.

\begin{Rem}
If $\sigma (t) = \sigma_1 t^{-2}$ for any $t \geq r_0$ with $\sigma_1 \in [0,1/4)$, then we have two fundamental solutions $t^{\lambda}$ and $t^{1- \lambda}$ of $y''(t) + \sigma (t) y(t) = 0$ with $\lambda = (1-\sqrt{1-4 \sigma_1})/2$. Then $\zeta _j (t)$ can be written as a linear combination of $t^{\lambda}$ and $t^{1- \lambda}$. In various places, we use the inequality $ \left| e^{i \zeta _1(t)|x|^2/(\zeta _2(t)) } -1 \right| \leq C |x|^{2 \theta} (\zeta _1/\zeta_2)^{\theta}$ for $\theta \in [0,1]$ and deduce the decay in $t$ from this term. To obtain the decay, it is necessary that $\zeta_1/\zeta_2$ decays in $t$ and hence there are some additional conditions in \eqref{K13}.
\end{Rem}

In order to present the main result, let us briefly recall the decomposition of the nonlinearity used by \cite{MM2, MMU}.
We identify a homogeneous nonlinearity $F$ and $2\pi$-periodic function $g$ as follows:
A homogeneous nonlinearity $F$ is written as
\begin{equation}\label{eq:id1}
	F(u) = |u|^{1 + \frac{2}{d(1-\lambda)}} F\( \frac{u}{|u|} \).
\end{equation}
We then introduce a $2\pi$-periodic function $g(\theta)=g_F(\theta)$ by $g_F(\theta) = F(e^{i\theta})$.
Conversely, for a given $2\pi$-periodic function $g$, one can construct a homogeneous nonlinearity $F=F_g:\C \to \C$ by
$F_g(u) = |u|^{\frac{2}{d(1-\lambda)}} g\( \arg u \)$ if $u\neq 0$ and
$F_g(u) = 0$ if $u=0$.
Since $g(\theta)$ is $2\pi$-periodic function, at least formally, $g(\theta)=\sum_{n\in\Z} {g}_n e^{in\theta}$ holds with \eqref{eq:gn}.
Thus, the expansion gives us \eqref{eq:exp} with $p_c = 2/d(1-\lambda)$.

\subsection*{Notations}
We introduce some notations used throughout this paper.
For any $p \geq 1$, $L^p = L^{p}(\R^d)$ denotes the usual Lebesgue space on $\R^d$.
Set $\J{a}=(1+|a|^2)^{1/2}$ for $a \in \C$ or $a\in \R^d$. Let $s$, $m \in \R$. The weighted Sobolev space and the homogeneous Sobolev space
on $\R^{d}$ are defined by $H^{m,s} = H^{m,s}(\R^d) = \{u \in \mathscr{S}' \mid \J{x}^s \J{i \nabla}^m u \in L^2 \}$ and 
$\dot{H}^s = \dot{H}^{s}(\R^d) = \{u \in \mathscr{S}' \mid |\nabla|^s u \in L^2\}$, respectively. 
Here $\mathscr{S}'$ is the space of tempered distributions.
We simply write $H^{m} = H^{m,0}$.
$\SCR{F}[u] = \widehat{u}$ is the usual Fourier transform of a function $u$ on $\R^d$ and $\SCR{F}^{-1}[u] = \check{u}$ is its inverse.  
$\norm{g}_{\mathrm{Lip}}$ stands for the Lipschitz norm of $g$.
We say a pair $(q,r)$ is an admissible if it satisfies 
\begin{align*}
\frac{1}{q} + \frac{d}{2r} = \frac{d}{4}, \quad q >2 , \quad r \geq 2. 
\end{align*} 

\subsection{Main results}

Throughout this paper, we always lie in the following condition for $\lambda = \lambda_d$: 
\begin{align}\label{K14}
\lambda _1 < 4- \sqrt{15}, \quad \lambda _2 < \frac15, \quad \lambda _3 < \frac{13- 2\sqrt{37}}{21}. 
\end{align}
Remark that it is expected that we can weaken the condition \eqref{K14}, because the case $\lambda \in [0,1/2]$ is discussed in \cite{KM21, K21} on the initial value problem with the gauge-invariant nonlinearity.
Compared with \cite{MM2, MMU}, 
we assume a little stronger assumption for $\{g_n\}$ as follows:
\begin{Assu}\label{A2}
Assume that the nonlinearity $F \colon \C \to \C$ is a homogeneous of degree $1+2/d(1- \lambda)$ such that $g_n$ defined by \eqref{eq:gn} satisfies $g_0 =0$, $g_1 \in \R$ and
\begin{align*}
\sum_{n \in {\bf Z} } |n|^{1+ a(\lambda) + \eta} |g_n| < \infty 
\end{align*}
for some $\eta >0$, where $a(\lambda) = a_d(\lambda)$ is a monotone strictly increasing function on the domain \eqref{K14} defined by
\[
	a_1(\lambda) = \frac{6\lambda -\lambda^2}{4(1-2\lambda)}, \quad a_d(\lambda) = \frac{3 d \lambda}{4(1-2\lambda)} \quad (d=2,\, 3).
\]
Note that $a(\lambda) \in (0, 1/4)$ if $d=1$, $a(\lambda) \in (0,1/2)$ if $d=2$, and $a(\lambda) \in \(0, (2\sqrt{37}-11)/2 \)$ if $d=3$. Here $(2\sqrt{37}-11)/2 \sim 0.097$. In particular, $g$ is Lipschitz continuous.
\end{Assu}
Our main result is the following: 
\begin{Thm}\label{T1}
Suppose that $F$ satisfies Assumption \ref{A2} for some $\eta>0$. 
Let $d \leq 3$. Set $\delta = \delta_d >0$ such that
\begin{align}\label{K15}
\begin{aligned}
	\frac{1+4\lambda-\lambda^2}{2(1-2\lambda)} <{}& \delta_1 < \min \(1, \frac12 + 2a(\lambda) + 2\eta \) \quad (d=1), \\
	\frac{d(\lambda+1)}{2(1-2\lambda)} <{}& \delta_d < \min \(2 , 1+p_c, \frac{d}2 + 2a(\lambda) + 2\eta \) \quad (d=2,\, 3).
\end{aligned}
\end{align} 
Take $\delta^{\prime}= 1$ if $d=1$ and $\delta^{\prime} = \delta$ if $d=2$, $3$.  
Then, under Assumption \ref{A1}, there exists $\ep_0 = \varepsilon_0(\norm{g}_{{\rm Lip}}) >0$ such that for any $u_+ \in  H^{0,\delta'} \cap \dot{H}^{-\delta}$ with $\| \widehat{u_+} \|_{\infty} < \ep_0$, there exist a $T \geq r_0$ and an unique solution $u \in C([T, \infty) \, ; \, L^2({\bf R}^d))$ to \eqref{eq1} which satisfies 
\begin{align}
\sup_{t \geq T} t^{b-\lambda} \|u(t) -u_p (t) \|_{2} < \infty \label{main:1}
\end{align}
for any $b \in \(2\lambda, \lambda + \delta (1-2\lambda)/2 \)$, where
\begin{align}
	u_p (t) = (it)^{-\frac{d}{2}}e^{i\frac{|x|^2 \zeta'_2(t)}{2 \zeta_2(t)}} \widehat{u_+}\( \frac{x}{t}\) \exp\(-i \frac{g_1}{c_+} \left|\widehat{u_+}\(\frac{x}{t}\) \right|^{\frac{2}{d(1-\lambda)}} \log t\), \; c_+ = |c_{2}^{(0)}|^{\frac{1}{1-\lambda}}.
	\label{asym:1}
\end{align}
Moreover, there exists $\varepsilon_1 >0$ such that the solution satisfies 
\begin{align}
	\sup_{t \geq T }  t^{b-2 \lambda} \(\int_t^{\infty} \J{s}^{-\lambda} \|u(s) -u_p (s)\|_{r}^q\, ds \)^{\frac1q} < \infty \label{main:2}
\end{align}
for any admissible pair $(q,r)$ with $q_d \leq q$ and all $b \in \(2\lambda, \lambda + \delta (1-2\lambda)/2 \)$, 
where $q_d = 4$ if $d=1$ and $q_d = 2/(1-2 \varepsilon_1)$ if $d=2$, $3$. 
\end{Thm}

\begin{Rem}
We need the condition \eqref{K14} to take $\delta$ satisfying \eqref{K15}. In fact, if $d=1$, then \eqref{K15} reaches to $\lambda^2 - 12 \lambda +3 >0$ which yields $\lambda < 6-\sqrt{33} \sim 0.255$. 
When $d=2$, \eqref{K15} implies $\lambda < 1/5$. In $d=3$, we also have $21 \lambda^2 -26 \lambda +1 >0$. This yields 
$\lambda < (13- 2\sqrt{37})/21 \sim 0.0397$. 
\eqref{K14} is weaker than that in \cite{Ka2}. Therefore, we improve the condition of  $\delta$ as in \eqref{K15},
compared with \cite{Ka2}. 
However, since the lower bound of $\delta$ is larger than $d/2$, we impose additional decay rate $a(\lambda)$ in Assumption \ref{A2} (see Lemma \ref{lem:11}).
Also, when $\lambda =0$, Assumption \ref{A2}, \eqref{K15} and the condition for $b$ coincide those in \cite{MM2, MMU}.  
The lower bound of $\delta$ determined by that of $b$ in Proposition \ref{propm:1}.
\end{Rem}

\begin{Rem}
$\sigma(t) \equiv 0$ satisfies Assumption \ref{A1}. 
In fact, we can take $\zeta_1 (t) = 1$, $\zeta_2(t) = t$ for any $t \geq 0$.
Hence, Theorem \ref{T1} includes the free equation \eqref{eq:f} such as \cite{MM2, MMU}.
Furthermore, when $d=2$, $3$, compare with \cite{MM2, MMU}, 
we improve the regularity condition of $u_+$, by applying Proposition \ref{pro:a2} and Proposition \ref{pro:a4}. 
\end{Rem}

\begin{Rem}
If we take $\varepsilon_1 = 0$, then the admissible pairs $(q_d, r_d)$, $d=2,3$ defined by \eqref{K2} below, correspond to the so-called end-point of admissible pairs. 
To relax the range of $\lambda$ more than that in \cite{Ka2}, we need to choose the admissible pair near the end-point.
Since the end-point Strichartz estimate for the propagator of $H_0(t)$ in $d=2,3$ has not been proven, for a technical reason, $\varepsilon_1$ is demanded.
$\varepsilon_1$ is precisely specified by \eqref{M2} and \eqref{M3} in Section \ref{sec:2}.
\end{Rem}

\begin{Rem}
When $F(u)$ satisfies Assumption \ref{A2} and $g_1 \neq 0$, Theorem \ref{T1} tells us that \eqref{eq1} admits a nontrivial solution that behaves like a free solution with the logarithmic phase correction. Namely, such a nonlinearity is long-range.  
The typical example of the nonlinearity is $F(u) = |\Re u|^{\frac{2}{d(1-\lambda)}} \Re u$. In fact, the corresponding periodic function is  
\[
	g(\theta) = |\cos \theta|^{\frac{2}{d(1-\lambda)}} \cos \theta, \quad g_n = O(|n|^{-\frac{2}{d(1-\lambda)}-2})
\]
as $n \to \I$.
Further, under Assumption \ref{A2} and $g_1 = 0$, Theorem \ref{T1} implies that \eqref{eq1} admits a nontrivial solution with the asymptotic profile 
\[
	(it)^{-\frac{d}{2}}e^{i\frac{|x|^2 \zeta'_2(t)}{2 \zeta_2(t)}} \widehat{u_+}\( \frac{x}{t}\)
\]
corresponding to that of the free solution. Thus, in this case, the nonlinearity is short-range. 
We then have the example $F(u) = |\Re u|^{\frac{2}{d(1-\lambda)}} \Re u - i |\Im u|^{\frac{2}{d(1-\lambda)}} \Im u$ with
\[
	g(\theta) = |\cos \theta|^{\frac{2}{d(1-\lambda)}} \cos \theta - i |\sin \theta|^{\frac{2}{d(1-\lambda)}} \sin \theta, \quad g_n = O(|n|^{-\frac{2}{d(1-\lambda)}-2})
\]
as $n \to \I$. For the computation of $g_n$, we refer to \cite{MMU}.
\end{Rem}

\begin{Rem} \label{rem:18}
Precisely, the uniqueness assertion of the solution is in the following sense: If $\widetilde{u} \in C([T, \infty)\, ;\, L^{2}(\R^d))$ solves \eqref{eq1} and satisfies \eqref{main:1} and \eqref{main:2} for some $\widetilde{b} > 0$ satisfying \eqref{Mb:1} in Proposition \ref{prop:31}, 
then $u = \widetilde{u}$.
\end{Rem}

\subsection{Strategy of the proof of main results and setting}

The strategy of the proof of the main results is based on that of \cite{MM2, MMU, HNW}.
We now let $U_0(t,s)$ be a propagator for $H_0 (t)$, that is, the family of unitary operators $\{ U_0(t,s)\}_{(t,s) \in {\bf R}^2} $ on $L^2({\bf R}^n)$ such that for all $t,s,\tau \in {\bf R}$,
\begin{align*}
& i \partial _t U_0(t,s) = H_0(t) U_0(t,s), \quad i \partial _s U_0(t,s) = -U_0(t,s) H_0(s), \\
& U_0(t, \tau) U_0(\tau, s) = U_0(t,s), \quad U_0(s,s) = \mathrm{Id}_{L^2({\bf R}^n)}, \quad U_0(t,s) \D{H_0(s)} \subset \D{H_0(s) }
\end{align*}
hold on $\D{H_0(s)}$. 
By using $\zeta_1 (t)$ and $\zeta_2 (t)$, the following factorization formula of $U(t, 0)$ can be originally obtained by Korotyaev \cite{Ko}, and 
Kawamoto and Muramatsu \cite{KM21} rewrites the formula to employ the factorization technique developed by \cite{HN98} (see also \cite{Ca1, Ca2, Ka2}).

\begin{Lem}[\cite{Ko, KM21}] \label{mdfm:1}
For $\phi \in \SCR{S}({\bf R}^d)$, let us define
\begin{align*}
\left( \CAL{M}(\tau) \phi \right) (x) = e^{\frac{i|x|^2}{2 \tau}} \phi (x), \quad
\left(
\CAL{D}(\tau) \phi
\right) (x) = \frac{1}{(i \tau)^{n/2}} \phi (x/ \tau).
\end{align*}
Then the following holds:
\begin{align}\label{mdfm}
U_0(t,0) = \CAL{M} \left(  \frac{\zeta _2(t)}{\zeta _2 '(t)} \right) \CAL{D} (\zeta _2 (t)) \SCR{F} \CAL{M} \left( \frac{\zeta _2 (t)}{ \zeta _1 (t)} \right).
\end{align}
\end{Lem}

Hereafter we use the notation 
\begin{align*}
\CAL{M}_1 (t) = \CAL{M} \left(  \frac{\zeta _2(t)}{\zeta _2 '(t)} \right), \quad \CAL{M}_2 (t) =\CAL{M} \left( \frac{\zeta _2 (t)}{ \zeta _1 (t)} \right).
\end{align*}
Also, set 
\begin{align*}
u_p (t) ={}& \CAL{M}_1(t) \CAL{D}(\zeta _2 (t)) \widehat{w} (t), \quad
	\widehat{w} (t)  = \widehat{u_+} \exp \left( -i g_1 |\widehat{u_+}|^{p_c} \log t/c_+ \right),
\end{align*}
and
\begin{align*}
\CAL{G}(u(t)) ={}& g_1 |u(t)|^{p_c} u(t), \quad 
\CAL{N}(u(t)) = \sum_{n \neq 0,1} g_n |u(t)|^{1 + p_c -n}u(t)^n, \\
R(t) ={}& \CAL{M}_1 (t) \CAL{D}(\zeta _2 (t)) \left( \SCR{F}\CAL{M}_2 (t) \SCR{F}^{-1} -1 \right).
\end{align*}
By using the formulation obtained by Kawamoto \cite{Ka2} (see also \cite{HNW}),
we rewrite \eqref{eq1} as the integral equation 
\begin{align} \label{inteq:1}
u(t) - u_p(t) = i \int_t^{\infty} U_0(t,s) \left( 
F(u(s)) - F(u_p(s))
\right)  ds + \CAL{A}(t) + \CAL{E}_{\mathrm{r}} (t) + \CAL{E}_{\mathrm{nr}} (t), 
\end{align}
where 
\begin{align*}
\CAL{A}(t) ={}& i \int_t^{\infty} U_0(t,0) \SCR{F}^{-1} \left( \frac{c_+s}{\zeta _2 (s)^{1/(1- \lambda)} } -1 \right) \mathcal{G}(\widehat{w}(s)) \frac{ds}{c_+ s}, \\[5pt]
\CAL{E}_{\mathrm{r}} (t) ={}& R(t) \widehat{w}(t) -i \int_t^{\infty} U_0(t,s) R(s) \mathcal{G}(\widehat{w}(s)) \frac{ds}{\zeta_2(s)^{1/(1-\lambda)}}, 
\end{align*}
and 
\begin{align*}
\CAL{E}_{\mathrm{nr}} (t) &= i \int_t^{\infty} U_0(t,s) \CAL{N} (u_p (s))\, ds.
\end{align*}

In the followings,  for all $q \in [1,\infty]$, we use the notation $\| \cdot \|_{q}$ as the norm on $L^q({\bf R}^d)$ for short. 
Moreover, for all $q$, $r \in [1,\infty]$ and $a<b$, we define the time-weighted Bochner-Lebesgue space $L^q_{\lambda} ((a,b) ; L^r({\bf R}^d))$ (see Kawamoto-Yoneyama \cite{KY1}) by
\begin{align*}
L_{\lambda}^q ((a,b) ; L^r ({\bf R}^d)) = \left\{ 
f \in \SCR{S}' ({\bf R}^{1+d}) \, \middle| \, 
 \| f \|_{L^q_{\lambda} ((a,b) ; L^r({\bf R}^d))} < \infty
\right\},
\end{align*}
where 
\begin{align*}
\| f \|_{L^q ((a,b) ; L^r({\bf R}^d), \lambda )} = \left( 
\int_{a}^{b} \J{t}^{-\lambda} \left\| 
f(t, \cdot)
\right\|_{r}^q\, dt 
\right) ^{1/q}.
\end{align*} 
Here we denote a function space on which we can construct the contraction mapping. 
In \cite{HNW, MM2, MMU}, the space equips the same time-weight $t^b$. 
On the other hand, in our case, the different weight needs to be equipped in the function space associated with the Strichartz norm.
This fact was found by \cite{Ka2} and hence we also employ the same weighted space. 
Let us define
\begin{align*}
\left\| 
f
\right\|_{q, r, \lambda, \tau } = \left( \int_{\tau}^{\infty}
\J{s}^{-\lambda} \left\| f(s ,\cdot ) \right\|_{r}^{q} ds 
\right) ^{1/q}, \quad 
\left\| 
f
\right\|_{\infty, r, \lambda, \tau } = \sup_{s \geq \tau}
\J{s}^{-\lambda} \left\| f(s ,\cdot ) \right\|_{r}. 
\end{align*} 
We introduce a complete metric space 
\begin{align}
	\begin{aligned}
	X_{T,b,R} \coloneqq{}& \left\{ \phi \in C\left( [T, \infty) \, ; \, L^2({\bf R}^d) \right) \; ;\; \left\|\phi - u_p \right\|_{X_{T,b}} \leq R \right\}, \\
	\left\|\phi\right\|_{X_{T,b}} \coloneqq{}& \sup_{\tau \geq T } \tau^{b} \| \phi \|_{\infty, 2,\lambda, \tau} + \sup_{ \tau \geq T }  \tau ^{b-2 \lambda} \|\phi \|_{q, r, \lambda, \tau}, \\
	d(u,v) \coloneqq{}& \norm{u-v}_{X_{T,b}}
	\end{aligned}
	\label{K1}
\end{align}
for $R>0$ and $T \geq r_0$, where $(q, r) = (q_d, r_d)$ is a certain admissible pair defined by
\begin{align}
(q_1, r_1) = (4, \infty), \quad (q_2,r_2) = \left(\frac{2}{1-2 \varepsilon_1}, \frac{1}{\varepsilon_1} \right), 
\quad (q_3,r_3) = \(\frac{2}{1-2 \varepsilon_1}, \frac{6}{1+4 \varepsilon_1}\).
	\label{K2}
\end{align}

The main step of the proof of Theorem \ref{T1} is to show that the right-hand side of \eqref{inteq:1} is negligible for large time.
The first term does not matter because it can be controlled by the Strichartz estimate associated with the propagator $U_0(t,s)$ stated in Section \ref{sec:2}.
However, the restriction of $\lambda$ comes from this term, due to a use of the Strichartz estimate.
The terms $\CAL{A}(t)$ and $\CAL{E}_{\mathrm{r}} (t)$ are also harmless as in previous works \cite{Ka2, HNW}. 
The main issue of the proof is to estimate the non-resonant term as follows:

\begin{Prop} \label{nonres}
Let $u_+ \in H^{0, \delta'} \cap \dot{H}^{-\delta}$.
Then there exists $C=C(\norm{u_+}_{H^{0, \delta'} \cap \dot{H}^{-\delta}}) > 0$ such that
\begin{align*}
	\norm{\mathcal{E}_\mathrm{nr}}_{\I, 2, \lambda, \tau} + \tau^{-2\lambda} \norm{\mathcal{E}_\mathrm{nr}}_{q, r, \lambda, \tau} 
	\leq C\J{\tau}^{-\frac{\delta}2 (1-2\lambda) - \lambda} \J{g_1 \log \tau}^{\lceil \delta \rceil} \sum_{n \in \Z} |n|^{1+a(\lambda) + \eta} |g_n|
\end{align*}
for all $\tau \geq r_0$. 
\end{Prop}

To establish Proposition \ref{nonres}, we follow the argument of \cite{MM2, MMU}. 
However, one needs to adapt a regularizing operator and a factorization identity for the propagator, developed by \cite{HNW}, to \eqref{eq1}. Then, additional terms appear due to the factorization formula of the propagator, compared with \cite{MM2, MMU}.
Further, to get the necessary time decay for the additional terms, we need to take it into account for a regularizing operator to degenerate in the low-frequency region.

The rest of the paper is organized as follows:
In section \ref{sec:2}, we introduce the Strichartz estimate with time weight for the propagator and prove the existence of solutions near the given asymptotic profile corresponding to the main theorem in an abstract form.
Then Section \ref{sec:3} is devoted to summarize useful estimates.
In the section, we present estimates for the fractional derivative of $\widehat{w}$ necessary to improve the regularity condition for the final data in $d=2$, $3$. 
One also states the estimates for $\CAL{A}(t)$ and $\CAL{E}_{\mathrm{r}} (t)$ given by \cite{Ka2, HNW}, and
adapts the regularizing operator and the factorization identity for the propagator to \eqref{eq1}.
Furthermore, Proposition \ref{nonres} are shown and our main result is established in Section \ref{sec:nre}.
Finally, Appendix provides us to prove the estimates for $\widehat{w}$ given in Section \ref{sec:3}. 
Remark that the estimates for a general form of $\widehat{w}$ are treated in the section.


\section{Existence of solutions near the asymptotic profile} \label{sec:2}

Let us introduce the so-called Strichartz estimate for associated propagator $U_0(t,s)$ which was firstly obtained by Kawamoto and Yoneyama \cite{KY1} for the case of the restricted coefficients $\sigma (t)$. Subsequently, Kawamoto \cite{Ka4} handles the Strichartz estimate under the more generalized condition including the assumption \ref{A1}.
\begin{Lem}[\cite{KY1, Ka4}] \label{Str:1}
Let $(q,r)$ and $(\tilde{q}, \tilde{r} )$ be admissible pairs and $s'$ denotes the H\"{o}lder conjugate of $s$, i.e., $1/s + 1/s' =1$.
Then for $\phi \in L^2({\bf R}^d)$ and $F \in L_{\lambda} ^{\tilde{q}'}((a,b) ; L^{\tilde{r}'}({\bf R}^d), - \lambda)$, there exists $C>0$ which is independent of $a,b$ such that it holds that
\begin{align*}
\left\| 
U_0(t,0) \phi
\right\|_{L^q((a,b);L^r({\bf R}^d) , \lambda)} \leq C \| \phi \|_2
\end{align*}
and 
\begin{align*}
\left\| 
\int^t U_0(t,s) F(s) ds 
\right\|_{L^q ((a,b);L^r({\bf R}^d), \lambda)} \leq C \left\| F \right\|_{L^{\tilde{q}'}((a,b);L^{\tilde{r}'}({\bf R}^d) , - \lambda)}.
\end{align*}
\end{Lem}

In $d=2$, $3$, we fix $\varepsilon_1>0$ is a sufficiently small constant with 
\begin{align}
	\frac{d(\lambda+1)}{2(1-2\lambda)} + \frac{d (1-\lambda) \lambda \varepsilon_1}{1-2\lambda} < \delta
	\; \Leftrightarrow \;  \frac{d(\lambda+1)}{4} + \frac{d (1-\lambda) \lambda}{2} \varepsilon_1 < \frac{\delta}{2} (1-2\lambda),
	\label{M2}
\end{align}
and
\begin{align}
	\lambda+1 + 2 (1-\lambda) \lambda \varepsilon_1
	>{}& \frac{2(1-\lambda)}{d+(1-2 d) \lambda- d(1-\lambda) \varepsilon_1}
	\label{M3}
\end{align}
When $\lambda$ satisfies \eqref{K14}, we have
\[
	\frac{d(1-\lambda)}{2d+2(1-2 d) \lambda} < \frac{d}{4}(\lambda+1),
\]
which ensures that one can take $\varepsilon_1$ satisfying \eqref{M3}.

In what follows, we assume $\lambda >0$, since the case $\lambda =0$ can be proven in a similar way to \cite{MM2, MMU}.
Let us consider the integral equation
\begin{align} \label{inteq:2}
u(t) - u_p(t) = i \int_t^{\infty} U_0(t,s) \left( F(u(s)) - F(u_p(s)) \right)  ds + \mathcal{E}(t),
\end{align}
where $u_p$ is the given asymptotic profile \eqref{asym:1} and $\mathcal{E}(t)$ is an external term.
Remark that our equation \eqref{inteq:1} is of the form \eqref{inteq:2}.
As for the existence of solution to \eqref{inteq:2}, we have the following:

\begin{Prop} \label{prop:31}
Suppose that $g$ is Lipschitz continuous. Let $(q,r)$ be the certain admissible pair as in \eqref{K2}.
Set $\widehat{u_+} \in L^{\infty}$ and $u_p$ be as in \eqref{asym:1}.
Then there exists a constant $\varepsilon_0 = \varepsilon_0(\norm{g}_{\mathrm{Lip}})>0$ such that
if $\norm{\widehat{u_+}}_{L^\I} \leq \varepsilon_0$
and if an external term $\mathcal{E}$ satisfies
$\norm{\mathcal{E}}_{X_{T_0,b}} \le M$ for some $T_0 \geq 1$, $M>0$, and $b = b_{d}$ with
\begin{align}
\begin{aligned}
	\frac{1+8 \lambda - \lambda^2}{4} <{}& b_1 < \lambda + \frac{\delta}{2}(1-2\lambda), \\
	\frac{d}{4}(\lambda+1)+\lambda + \frac{d}{2} (1-\lambda) \lambda \varepsilon_1 <{}& b_d < \lambda + \frac{\delta}{2}(1- 2 \lambda) 
	\quad (d=2,\, 3),
\end{aligned}
\label{Mb:1}
\end{align}
then there exists a $T_1= T_1(M, \norm{g}_{\mathrm{Lip}}, b) \geq T_0$ such that \eqref{inteq:2} admits a unique solution $u(t)$ in $X_{T_1,b, 2M}$.
Moreover, for any admissible pair $( \widetilde{q}, \widetilde{r})$ with $q \leq \widetilde{q}$, and $\tilde{b} \leq b$, 
there exists $C_0 = C_0(M) >0$ such that 
the solution satisfies
\[
	\sup_{t \geq T_1} t^{ \widetilde{b}-2 \lambda} \norm{u-u_p}_{\widetilde{q}, \widetilde{r}, \lambda, t}
	\le C_0 + \sup_{t \geq T_1} t^{\widetilde{b}-2 \lambda} \norm{\mathcal{E}}_{\widetilde{q}, \widetilde{r}, \lambda, t}.
\] 
\end{Prop}

To prove Proposition \ref{prop:31},
we now estimate the $X_{T,b}$ norm of the first term of the right-hand side of \eqref{inteq:2}. This term has been estimated in Lemma 3.3 of \cite{Ka2}. However, the proof employs a very simple approach to avoid the complex calculations and then the strong restrictions for the higher and lower bound of $b$ are required. Hence the aim here is to relax these restrictions.
\begin{Prop} \label{propm:1}
Assume that $g$ is Lipschitz continuous. Then there exists $\widetilde{\varepsilon}_0 >0$ and $C>0$ such that 
\begin{multline*}
	\left\| \int_t^{\infty} U_0(t,s) \left( F(u(s)) -F(u_p(s)) \right) ds \right\|_{X_{T,b}} \\
	\leq C \norm{g}_{{\rm Lip}} \left\| u-u_p\right\|_{X_{T,b}} \( T^{-\widetilde{\varepsilon}_0} \left\| u-u_p\right\|_{X_{T,b}}^{p_c} + \left\| \widehat{u_+} \right\|_{\infty}^{p_c} \)
\end{multline*}
for any $b=b_d$ satisfying \eqref{Mb:1}.

\end{Prop}
\begin{Rem}
The upper bound of $b$ in \eqref{Mb:1} comes from \eqref{M2}, and is not necessary to prove Proposition \ref{propm:1}.
However, the bound is crucial to prove Theorem \ref{T1}.
\end{Rem}

\begin{proof} 
The nonlinear terms can be decomposed as $F(u) - F(u_p) = F^{(1)} (u) + F^{(2)}(u)$, where
\begin{align}
	F^{(1)}(u) = \chi_{\{|u_p| \leq |u -u_p|\}} \(F(u) - F(u_p)\), \quad
	F^{(2)}(u) = \chi_{\{|u_p| \geq |u -u_p|\}} \(F(u) - F(u_p)\),
	\label{pr:31c}
\end{align}
and $\chi_A$ is a characteristic function on $A \subset \R^{1+d}$. 
Since $g$ is Lipschitz continuous, 
we see from \cite[Appendix A]{MM2} that
\[
	|F(u) - F(u_p)| \leq C \norm{g}_{{\rm Lip}} \( |u-u_p|^{p_c+1} + |u_p|^{p_c}|u-u_p| \),
\]
and therefore 
\[
	|F^{(1)}(u)| \leq  C|u-u_p|^{p_c+1}, \quad |F^{(2)}(u)| \leq  C|u_p|^{p_c}|u-u_p|.
\]
By the same calculations in \cite{Ka2}, we obtain for $b < \lambda + \delta (1-2 \lambda) /2$, 
\begin{align*}
\left\| 
\int_t^{\infty} U_0(t,s)  F^{(2)} (u(s))\, ds 
\right\|_{\infty, 2, \lambda} \leq C \tau^{-b} \left\| \widehat{u_+} \right\|_{\infty}^{p_c} \left\| u-u_p \right\|_{X_T} 
\end{align*}
and 
\begin{align*}
\left\| 
\int_t^{\infty} U_0(t,s)  F^{(2)} (u(s))\, ds 
\right\|_{q, r, \lambda} \leq C \tau^{-(b - 2 \lambda)} \left\| \widehat{u_+} \right\|_{\infty}^{p_c} \left\| u-u_p \right\|_{X_T}. 
\end{align*}
Hence it is enough to show that there exists $\widetilde{\varepsilon}_0>0 $ such that
\begin{align}
\left\| 
\int_t^{\infty} U_0(t,s)  F^{(1)} (u(s))\, ds 
\right\|_{\infty, 2, \lambda} \leq C \tau^{- b  -\widetilde{\varepsilon}_0} \left\| u-u_p\right\|_{X_T}^{1 + p_c} \label{Ma:1}
\end{align}
and 
\begin{align}
\left\| 
\int_t^{\infty} U_0(t,s)  F^{(1)} (u(s))\, ds 
\right\|_{q,r, \lambda} \leq C \tau^{- (b-2 \lambda)  - \widetilde{\varepsilon}_0 } \left\| u-u_p\right\|_{X_T}^{1 + p_c} \label{Ma:2}
\end{align}
for any $\tau \geq r_0$.
By the Strichartz estimate, we estimate 
\begin{align*}
\left\| \int_t^{\infty} U_0(t,s)  F^{(1)} (u(s))\, ds \right\|_{\infty,2, \lambda} 
&{}\leq C \left( \int_{\tau}^{\infty} \J{s}^{\lambda} \left\| \left| 
(u-u_p) (s)
\right|^{1 + p_c} \right\|_{\kappa'}^{\rho'}  ds \right) ^{1/ \rho'}
\\ 
&{}\eqqcolon C L(\tau)
 \end{align*}
with an admissible pair $(\rho,\kappa)$. Let $\theta _1 + \theta _2 = 1 + p_c$, and then by the H\"{o}lder inequality 
\begin{align*}
\left\| 
\left| 
u-u_p
\right|^{1 + p_c}
\right\|_{\kappa'} \leq \left\|u-u_p \right\|_{r}^{\theta _1} \left\| u-u_p \right\|_{2}^{\theta _2}
\end{align*}
with
\begin{align*}
\frac{1}{\kappa'} = \frac{\theta_1}{r} + \frac{\theta _2}{2} \; \Leftrightarrow \; \frac{2}{\rho} + \frac{2}{q} \cdot \theta_1 = \frac{1}{1-\lambda}, \quad \theta _2 = 1+p_c -\theta_1,
\end{align*}
where $\theta_1 \in [\theta_0, p_c+1]$ with $\theta_0 = \frac{1+\lambda}{1-\lambda} \in (1, 2)$ if $d=1$, otherwise, $\theta_1 \in (\theta_0, p_c+1]$ with $\theta_0 = \frac{\lambda}{1-\lambda}$, in view of the conditions $0 \leq 1/\kappa' -1/2 = 2/d\rho  \leq 1/2$ for $d=1$ and $0 \leq 1/\kappa' -1/2 = 2/d\rho  < 1/d$ for $d=2,3$.
This yields the restrictions 
\begin{align}\label{K4}
\frac{1}{r} \leq \frac{\theta_1 - \theta_0}{2\theta _1}, \quad \frac{1}{q} \geq \frac{\theta_0}{4\theta _1}
\end{align}
for $d=1$ and 
\begin{align} \label{K5}
\frac{1}{r} < \frac{d \theta_1 - 2 \theta_0}{2d \theta _1}, \quad \frac{1}{q} > \frac{\theta_0}{2 \theta _1}. 
\end{align}
for $d=2,3$.
Hence by the H\"older inequality, we have 
\begin{align*}
L(\tau) \leq{}& C \left( \int_{\tau}^{\infty} \J{s}^{\lambda} \left\| u-u_p \right\|_{r}^{\theta _1 \rho '} \left\| u-u_p \right\|_{2}^{\theta _2 \rho '} ds  \right)^{1/ \rho'}
\\ 
\leq{}& C \left( \sup_{s \geq \tau} \J{s}^{- \lambda + b} \left\| u-u_p \right\|_2 \right)^{\theta _2} \left\| \J{s}^{\lambda / \rho ' -b \theta _2 + \lambda \theta _2 } \left\| u-u_p \right\|_{r}^{\theta _1} \right\|_{L^{\rho '}([ \tau, \infty))}. 
\end{align*}
Using the H\"{o}lder inequality again, one has 
\begin{align*}
L(\tau) \leq C \left\| u-u_p \right\|_{X_T}^{\theta _2} \left\| \J{s}^{\lambda / \rho' -b \theta _2 + \lambda \theta _2 + \lambda \theta _3 / \rho '} \right\|_{L^{P_1}([\tau, \infty))} \left\| \J{s}^{- \theta _3 \lambda/ \rho'} \left\| u-u_p \right\|_{r}^{\theta _1}  \right\|_{L^{P_2}([\tau, \infty))},
\end{align*}
where $1/P_1+1/P_2 = 1/\rho'$, $\theta_1 P_2 = q$ and $\theta _3 q = \theta _1 \rho '$.
Hence this implies 
\begin{align}\label{K3}
L(\tau) \leq C  \left\| u-u_p \right\|_{X_T}^{\theta _2} \left\| \J{s}^{\lambda / \rho' -b \theta _2 + \lambda \theta _2 + \lambda \theta _1 / q} \right\|_{L^{P_1}([\tau, \infty))}  \left\| u-u_p \right\|_{q, r, \lambda}^{\theta _1}. 
\end{align}
In order to guarantee the right-hand side of \eqref{K3} is bounded, it is necessary to assume 
\begin{align*}
\frac{\lambda}{\rho'} -b \theta _2 + \lambda \theta _2 + \frac{\lambda \theta_1}{q} + \frac{1}{P_1}   < 0.
\end{align*}
Since
\begin{align*}
\frac{1}{\rho'} = 1 + \theta _1 \left( \frac{d}{4} - \frac{d}{2r} \right) + \frac{d-d (\theta _1 + \theta_2)}{4} = 1 + \frac{\theta _1}{q} - \frac{d p_c}{4}, 
\end{align*}
we have 
\begin{align*}
\frac{\lambda}{ \rho'} + \frac{\lambda \theta _1}{q} + \frac{1}{P_1} &= \frac{\lambda +1}{\rho '} + \frac{\lambda -1}{q} \theta _1
\\ &= \lambda +1 + \frac{2 \lambda \theta _1}{q} - \frac{d p_c (1 + \lambda)}{4} \eqqcolon \Lambda_d, 
\end{align*}
where we remark that by using \eqref{K4} and \eqref{K5}, one finds 
\begin{align*}
\Lambda _d &\geq \lambda +1 - \frac{p_c( 1+ \lambda )}{4} + \frac{\lambda (p_c-1)}{2} = \frac{1 + \lambda}{2}
\end{align*} 
for $d=1$ and 
\begin{align*}
\Lambda _d &> \lambda + 1  - \frac{d p_c (1 + \lambda)}{4} + \frac{d(1 + p_c) -2 -d}{2} \lambda = \frac12 
\end{align*}
for $d=2,3$. Hence, \eqref{K3} holds if $b$ satisfies  
\begin{align}
b >  \frac{\Lambda _d}{ \theta _2} + \lambda. \label{Mc1}
\end{align} 
Under this condition, it follows that 
\begin{align*}
L(\tau) &\leq C \tau^{- b \theta _2 + \lambda \theta _2 +(\lambda +1) + 2 \lambda \theta _1/ q -d p_c (1+ \lambda) /4 + (-b + 2 \lambda ) \theta _1} \| u-u_p \|_{X_T}^{\theta _1+ \theta _2} 
 \\ & \leq C \tau^{- b - p_c b + \lambda \theta _2 +(\lambda +1) + 2 \lambda \theta _1/ q -d p_c (1+ \lambda) /4 +  2 \lambda  \theta _1} \| u-u_p \|_{X_T}^{1 + p_c}.
\end{align*}
Then there exists $\widetilde{\varepsilon}_0>0$ such that
\begin{align} \label{K9}
\left\| \int_t^{\infty} U_0(t,s)  F^{(1)} (u(s))\, ds \right\|_{\infty,2, \lambda} \leq  C \tau^{- b - \widetilde{\varepsilon}_0} \| u-u_p \|_{X_T}^{1 + p_c} .
\end{align}
for any $\tau	\geq r_0$, as long as
\begin{align*}
- p_c b + \lambda \theta _2 +(\lambda +1) + 2 \lambda \theta _1/ q -d p_c (1+ \lambda) /4 +  2 \lambda  \theta _1 < 0, 
\end{align*}
which is equivalent to
\begin{align}
b > \frac{\Lambda _d}{ p_c} + \frac{\lambda}{p_c} \left( \theta _2 + 2 \theta _1 \right). \label{Mc2}
\end{align}
Similarly to the above, noting $b > b -2\lambda$, there exists $\widetilde{\varepsilon}_0>0$ such that 
\begin{align}
	\norm{\int_t^{\I} U_0(t,s) F^{(1)}(u(s))\, ds}_{q,r,\lambda, \tau} \leq{}& C \tau^{-(b - 2\lambda )- \widetilde{\varepsilon}_0} \norm{u-u_p}_{X_T}^{p_c+1} \label{K9a}
\end{align}
for any $\tau \geq r_0$, as long as \eqref{Mc1} and \eqref{Mc2} hold. Thus we conclude \eqref{Ma:1} and \eqref{Ma:2}.

From now on, let us specify $\theta_1$. 
Compared with \cite{Ka2}, the suitable choice of $\theta_1$ enables us to find a weaker lower bound for $b$. Note that both $1/\theta_2 = 1/(1+p_c-\theta _1)$ and $\theta _2 + 2 \theta _1 = 1+p_c + \theta _1$ are monotone increasing with respect to $\theta _1$, and hence right-hand side of \eqref{Mc1} and \eqref{Mc2} take minimum on the smallest value of $\theta _1$. 
In $d=1$, we can take $\theta_1 = \theta_0$ of the minimum, because it is allowed us to take the end-point $(q, r) = (4, \infty)$.
However, in the case of $d=2$, $3$, since the end-points of the admissible pair are excluded, the minimum of $\theta_1$ can not be attained. Hence, we need to take $\theta_1 = \theta_0 + \varepsilon_1$ for the small $\varepsilon_1>0$.

\medskip

\noindent
\underline{\bf Case $d=1$.} \quad
We take $\theta_1 = p_c -1 = \theta_0$ and $\theta_2 = 2$. Then it enable us to take $q=4$ only.
In terms of \eqref{Mc1} and \eqref{Mc2}, we also notice that 
\begin{align*}
	\frac{\Lambda _d}{ \theta _2} + \lambda = \frac{1+5\lambda}{4}
	< \frac{1+8 \lambda - \lambda^2}{4} 
	=\frac{\Lambda _d}{ p_c} + \frac{\lambda}{p_c} \left( \theta _2 + 2 \theta _1 \right) 
\end{align*}
for any $\lambda$ satisfying \eqref{K14}. Hence \eqref{K9} and \eqref{K9a} are valid provided $b > (1+8 \lambda - \lambda^2)/4$.

\medskip
\noindent
\underline{\bf Case $d=2$, $3$.} \quad 
We take $\theta_{1} = \theta_0 + \varepsilon_1$ and $\theta_{2} = 1+p_{c} - \theta_0 - \varepsilon_1$ with the small $\varepsilon_1>0$ on \eqref{M2} and \eqref{M3}. 
Then one can choose $\( 1/q, 1/r \) = \( 1/2 - \varepsilon_1, (d - 2 + 4 \varepsilon_1)/2d \)$ satisfying \eqref{K5}. 
Hence \eqref{Mc1} gives us
\begin{align*}
	b > \lambda+\frac{d(1-\lambda)}{2\left\{ d+(1-2 d) \lambda- d(1-\lambda) \varepsilon_1 \right\}}.
\end{align*}
It then follows from \eqref{Mc2} that
\begin{align}
	b > \frac{d}{4}(\lambda+1)+\lambda+\frac{d}{2} (1-\lambda) \lambda \varepsilon_1. 
	\label{Mc3}
\end{align}
By \eqref{M3}, we conclude \eqref{K9} and \eqref{K9a} if $b$ satisfies \eqref{Mc3}, as desired.
\end{proof}

\begin{proof}[Proof of Proposition \ref{prop:31}]

From the integral equation \eqref{inteq:2}, we define the map
\begin{align*}
	\Psi(v(t)) - u_p(t) = i \int_t^{\infty} U_0(t,s) \left( F(v(s)) - F(u_p(s)) \right)  ds + \mathcal{E}(t)
\end{align*}
on $X_{T,b, R}$. 
Let us show that $\Psi \colon X_{T, b, R} \to X_{T, b, R}$ is a contraction map for some $R>0$.
Firstly, combining Proposition \ref{propm:1} with the assumption, one has
\begin{align}
	\norm{\Psi(v) - u_p}_{X_{T,b}} \leq C_1 \norm{g}_{{\rm Lip}} R \(T^{- \widetilde{\varepsilon}_0} R^{p_c} + \varepsilon_0^{p_c} \) + M \label{pr:31a}
\end{align}
for any $v \in X_{T, b, R}$ with $T \geq T_0$ and $R>0$. 
Let us show
\begin{align}
	d\(\Psi(v_{1}), \Psi(v_{2})\) \leq C_2 \norm{g}_{{\rm Lip}} \(T^{-\widetilde{\varepsilon}_0} R^{p_c} + \varepsilon_0^{p_c} \) d(v_{1}, v_{2}) \label{pr:31b}
\end{align}
for any $v_1$, $v_2 \in X_{T,b,R}$. We estimate
\begin{align*}
	|F(v_1) - F(v_2)| \leq C \norm{g}_{{\rm Lip}} |v_1-v_2| \( |v_1 - u_p|^{p_c} + |v_2 - u_p|^{p_c} \) + C \norm{g}_{{\rm Lip}} |u_p|^{p_c} |v_1 -v_2|.
\end{align*}
Hence we can decompose $F(v_1) - F(v_2)$ into two parts as in \eqref{pr:31c}. Similarly to the proof of Proposition \ref{propm:1}, using  
\[
	 	\norm{|v_1-v_2| |v_j -u_p |^{p_c}}_{\kappa'} \leq \left\|v_j-u_p \right\|_{r}^{\theta_1} 
	 	\left\| v_j-u_p \right\|_{2}^{\theta_2-1} \left\| v_1 - v_2 \right\|_{2}
\]
for $j=1$, $2$, one obtains \eqref{pr:31b}, 
where $\theta_1$ and $\theta_2$ are the same in the proof of Proposition \ref{propm:1}.
Set $R = 2M$. Taking $\varepsilon_{0} = \varepsilon\( \norm{g}_{{\rm Lip}} \)$ small enough such that
\[
	C_1 \norm{g}_{{\rm Lip}} \varepsilon_0^{p_c} \leq \frac{1}{4}, \quad 
	C_2 \norm{g}_{{\rm Lip}} \varepsilon_0^{p_c} \leq \frac{1}{4},
\]
we have 
\begin{align*}
	\norm{\Psi(v) - u_p}_{X_{T,b}} \leq{}& 2M, 
	\quad d\(\Psi(v_{1}), \Psi(v_{2})\) \leq \frac{1}{2} d(v_{1}, v_{2}), 
\end{align*}
as long as we choose $T = T_{1} \geq T_{0}$ satisfying $T_{1}^{-b} (2M)^{p_{c}} \leq \varepsilon_{0}^{p_{c}}$.
This implies that $\Psi \colon X_{T_1, b, 2M} \to X_{T_1, b, 2M}$ is a contraction map. 
Thus, we have a unique solution in $X_{T_{1}, b, 2M}$.

Let us show the latter assertion.
By the interpolation and Proposition \ref{propm:1}, we obtain
\begin{align*}
	&{}\norm{u -u_p}_{\widetilde{q}, \widetilde{r}, \lambda} \\
	\leq{}& \norm{\int_t^{\infty} U_0(t,s) \left( F(v(s)) - F(u_p(s)) \right)  ds}_{\widetilde{q}, \widetilde{r}, \lambda}^{1- a}
	\norm{\int_t^{\infty} U_0(t,s) \left( F(v(s)) - F(u_p(s)) \right)  ds}_{\infty, 2, \lambda}^{a} \\
	&{}+ \norm{\mathcal{E}}_{\widetilde{q}, \widetilde{r}, \lambda} \\
	\leq{}& C \J{t}^{-b+2\lambda} (2M) + \norm{\mathcal{E}}_{\widetilde{q}, \widetilde{r}, \lambda}
\end{align*}
for any $t \geq T_1$, where $a \in (0,1)$ with $\widetilde{q}(1-a) = q$. This implies
\[
	t^{\widetilde{b}-2\lambda} \norm{u -u_p}_{\widetilde{q}, \widetilde{r}, \lambda, t} 
	\leq C \J{t}^{\widetilde{b}-b} (2M) + \norm{\mathcal{E}}_{\widetilde{q}, \widetilde{r}, \lambda, t},
\] 
and thus the desired estimate holds.

\end{proof}


\section{Preliminary estimates} \label{sec:3}
In this section, we summarize some important estimates. In what follows, we assume that $\lambda = \lambda _d >0$ satisfies \eqref{K14}.
Define $\delta >0$ and $\delta ' >0$ as constants satisfying \eqref{K15}
and $\delta ' = 1$ for $d=1$, $\delta '= \delta$ for $d=2$, $3$.
Under \eqref{K14}, we find that $1+ p_c >2$ if $d=1$, $2$ and $1+ p_c < 2$ if $d=3$.
Moreover, \eqref{K15} always fails for $d \geq 4$ because $\lambda < 1/2$ and hence we only consider the case $d \leq 3$. 
Also, fix $n \in {\bf Z}$ and we handle estimates for each $n$

\subsection{Estimates for the fractional chain rule}

We denote the smallest integer $n_0$ such that $n_0 \geq \delta$ by $\lceil \delta \rceil$ for any $\delta \in \R$.
Firstly, the following is valid.
\begin{Prop} \label{pro:a1}
Suppose $u_+ \in H^{0, \delta}$. Then there exists $C >0$ such that
\begin{align*}
\norm{\widehat{w}(s)}_{H^{\delta}} \leq{}& C \J{g_1 \log s}^{\lceil \delta \rceil} \norm{u_+}_{H^{0, \delta}}\( 1 + \norm{\widehat{u_+}}_{\infty}^{p_c} \)^{\lceil \delta \rceil}
\end{align*}
for any $s \geq 1$.
\end{Prop}

In \cite{MM2, MMU}, in order to weaken the condition on the nonlinearity,  
they use an interpolation technique.
Roughly speaking, this technique enables us to improve the upper bound of $\norm{|\widehat{w}|^{1+p_c-n} \widehat{w}^n}_{H^{\delta}}$ from $O(n^{\lceil \delta \rceil})$ into $O(n^{\delta})$. 
However, they need stronger regularity conditions for the data in the technique.
We here improve the regularity condition by combining the fractional Leibniz rule and intermediate use of the interpolation.
The key of the removability is $\delta > 1$ so that it is only applied in $d=2$, $3$. 

\begin{Prop} \label{pro:a2}
Suppose $u_{+} \in H^{0, \delta^{\prime}}$. 
Take $\gamma = \delta$ if $d= 1$, $2$ and $\gamma \in \left( \delta, 1+ p_c \right)$ if $d=3$.
Then there exists $C >0$ such that
\begin{align*}
	\norm{|\widehat{w}(s)|^{p_c+1-n} \widehat{w}(s)^n}_{H^{\delta}} 
	\leq{}& C \left\langle n \right\rangle^{\gamma} \J{g_1 \log s}^{\lceil \delta \rceil} \norm{u_{+}}_{H^{0, \delta'}} \left\lVert \widehat{u_{+}} \right\rVert_{\infty}^{p_c} \( 1 + \norm{\widehat{u_{+}}}_{\infty}^{p_c} \)^{\lceil \delta \rceil}
\end{align*}
for any $s \geq 1$.
\end{Prop}

\begin{Prop} \label{pro:a4}
Suppose $u_{+} \in H^{0, \delta^{\prime}}$. 
Then there exists $C >0$ such that
\begin{align*}
	&{}\norm{\partial_t \left( |\widehat{w}(s)|^{p_c+1-n} \widehat{w}(s)^n  \right)}_{H^{\delta}} \\
	\leq{}& C |g_1| \left\langle n \right\rangle^{1 + \delta} s^{-1} \J{g_1 \log s}^{\lceil \delta \rceil} \norm{u_{+}}_{H^{0, \delta'}} \left\lVert \widehat{u_{+}} \right\rVert_{\infty}^{2p_c} \( 1 + \norm{\widehat{u_{+}}}_{\infty}^{p_c} \)^{\lceil \delta \rceil}
\end{align*}
for any $s \geq 1$.
\end{Prop}
We postpone the proof of the above propositions and show the generalized version of the propositions in the appendix.
A use of Proposition \ref{pro:a1} and Proposition \ref{pro:a2} gives us the time decay estimate of the resonant part $\mathcal{E}_r$.
The following tells us that the resonant part $\mathcal{E}_r$ is harmless:

\begin{Lem} \label{res:1}
Let $(q,r)$ be an admissible pair.
The followings hold:
\begin{align*}
	&{}\|R(\tau) \widehat{w}(\tau)\|_{\infty, 2, \lambda, \tau} + \tau^{-2\lambda} \|R(\tau) \widehat{w}(\tau)\|_{q, r, \lambda, \tau} \\
	\leq{}& C \tau^{-\lambda-\frac{\delta}{2}(1-2 \lambda)} \J{g_1 \log \tau}^{\lceil \delta \rceil} \norm{u_+}_{H^{0, \delta}}\( 1 + \norm{\widehat{u_+}}_{\infty}^{p_c} \)^{\lceil \delta \rceil},
\end{align*}
and
\begin{align*}
	&{}\left\|\int_{t}^{\infty} U_{0}(t, s) R(s) \mathcal{G}(\widehat{w}(s)) \frac{d s}{\zeta_{2}(s)^{1 /(1-\lambda)}}\right\|_{\infty, 2, \lambda, \tau} \\
	&{}+ \tau^{-2\lambda} \left\|\int_{t}^{\infty} U_{0}(t, s) R(s) \mathcal{G}(\widehat{w}(s)) \frac{d s}{\zeta_{2}(s)^{1 /(1-\lambda)}}\right\|_{q, r, \lambda, \tau}\\
	\leq{}& C\left|g_{1}\right| \tau^{-\lambda-\frac{\delta}{2}(1-2 \lambda)} \J{g_1 \log \tau}^{\lceil \delta \rceil} \norm{u_{+}}_{H^{0, \delta'}} \left\lVert \widehat{u_{+}} \right\rVert_{\infty}^{p_c} \( 1 + \norm{\widehat{u_{+}}}_{\infty}^{p_c} \)^{\lceil \delta \rceil}
\end{align*}
for any $\tau \geq r_0$.
\end{Lem}
\begin{proof}
The strategy is the same as in \cite{Ka2, HNW}. Hence we omit the proof.
\end{proof}

Regarding $\CAL{A}(t)$, we have the following:
\begin{Lem}[{\cite[Lemma 3.4]{Ka2}}] \label{mca:1}
Let $(q,r)$ be an admissible pair.
Assume $u_+ \in L^2$.
Under Assumption \ref{A1}, there exists $C > 0$ such that
\begin{align*}
	\norm{\CAL{A}(\tau)}_{\infty, 2, \lambda, \tau} + \tau^{-2\lambda} \norm{\CAL{A}(\tau)}_{q, r, \lambda, \tau} \leq C\J{\tau}^{-\frac{\delta}2 (1-\lambda) - \lambda} \norm{u_+}_{2} \norm{\widehat{u_+}}_{\infty}^{p_c}
\end{align*}
for any $\tau \geq r_0$.
\end{Lem}

\subsection{Properties for the time-dependent regularizing operator}
Let $\psi \in \mathscr{S}(\R^d)$, the Schwartz space.
In order to obtain the time decay property of the non-resonant part $\mathcal{E}_{\mathrm{nr}}$,
Masaki, the second author and Uriya \cite{MMU} (cf. \cite{MM2, HNW}) introduced a regularizing operator depending on $n$ and $t$ defined by 
\[
	\psi\( \frac{-i\nabla}{|n| \sqrt{t}} \) =\F^{-1} 
	\psi \( \frac{\xi}{|n| \sqrt{t}} \) \F.
\]
In this paper, in order to suit \eqref{eq1}, we modify the operator as 
\begin{align}
	\mathcal{K}_\psi = \mathcal{K}_\psi(t,n) \coloneqq \CAL{\psi} \( \frac{-i \nabla}{|n| t^{\rho/2}} \), \quad \rho = 1-2\lambda.
	\label{kpsi:1}
\end{align}
Arguing as in the proof of \cite[Lemma 2.1]{MMU}, we have the following:
\begin{Lem} \label{mol:1.1}
Let $\psi \in \mathscr{S}(\R^d)$.  
Fix $s\in \R $ and $\theta \in [0,2]$.
Assume $\nabla \psi(0)=0$ if $\theta \in (1,2]$. For any $t>0$ and $n\neq 0$, the followings hold.
\begin{enumerate}
\renewcommand{\labelenumi}{(\roman{enumi})}
\item
$\mathcal{K}_\psi$ is a bounded linear operator on $L^2$ and satisfies $\norm{\mathcal{K}_\psi}_{\mathcal{L}(L^2)} \leq \norm{\psi}_{\I}$.
Further, $\mathcal{K}_\psi$ commutes with $\nabla$. In particular, $\mathcal{K}_\psi$ is a bounded linear operator on $\dot{H}^s$ and satisfies $\norm{\mathcal{K}_\psi}_{\mathcal{L}(\dot{H}^s)} \leq \norm{\psi}_{\I}$.
\item 
$\mathcal{K}_{\psi} - {\psi}(0)$ is a bounded linear operator from $\dot{H}^{s+\theta}$ to $\dot{H}^{s}$ with norm
$$\norm{\mathcal{K}_\psi-\psi(0)}_{\mathcal{L}(\dot{H}^{s+\theta} ,\dot{H}^{s})} \leq Ct^{-\frac{\rho}{2}\theta}|n|^{-{\theta}}.$$
\end{enumerate}
\end{Lem}


\subsection{Factorization of propagator}

We here show the factorization identity for the propagator necessary to analyze the non-resonant part of the nonlinearity.
\begin{Lem} \label{lem:41}
Let $ a ,b \in {\bf R}$ with $4ab \neq -1$. Then for all $\phi \in L^2({\bf R}^d)$,
\begin{align} \label{K16}
e^{ia \Delta } e^{ib |x| ^2} \phi = (i)^{d/2} e^{4ab i |x | ^2 / (1 + 4ab)} e^{a (1 + 4ab) i \Delta } \CAL{D} (1 + 4ab) \phi
\end{align}
holds.
\end{Lem}
\begin{proof}
We mimic the approach of \cite[Lemma A.2]{HNW}. For $\phi \in L^2({\bf R}^d)$, 
\begin{align*}
e^{ia \Delta } \left( e^{ib |x| ^2} \phi \right) (x) = \left(  \frac{1}{4 \pi i a} \right)^{d/2} \int_{{\bf R}^d } e^{i|x-y|^2/(4a)} e^{ib |y|^2} \phi (y) dy
\end{align*}
holds. Here 
\begin{align*}
\frac{|x-y|^2}{4a} + b |y|^2 &= \frac{|x|^2}{4a} - \frac{x \cdot y}{2a} + \frac{|y|^2}{4a} + b|y| ^2 
\\ &= 
\frac{|x| ^2}{4a} + \frac{1 + 4ab}{4a} \left( |y|^2 - \frac{2 x \cdot y}{1 + 4ab} \right) 
\\ &= 
\frac{4ab |x| ^2}{4a (1 + 4ab)} + \frac{1}{4a (1 + 4ab)} \left| x- (1 + 4ab) y \right|^2. 
\end{align*}
Since
\begin{align*}
& \int_{{\bf R}^d} e^{i \left| x-  (1 + 4ab)y \right|^2/(4a (1 + 4ab))} \phi (y) dy 
\\ 
&= \frac{1}{(1 + 4ab)^d} \int_{{\bf R}^d} e^{i \left| x- z \right|^2/(4a (1 + 4ab))} \phi \left( \frac{z}{1 + 4ab} \right) dz  
\\ 
&= 
\frac{1}{(-i( 1 + 4ab))^{d/2}} \int_{{\bf R}^d} e^{i \left| x- z \right|^2/(4a (1 + 4ab))} \left( \CAL{D} \left( 1 + 4ab \right) \phi \right) (z) dz
\\
&= 
\frac{(4 \pi i (a(1 + 4ab)))^{d/2}}{(-i( 1 + 4ab))^{d/2}} e^{a(1+4ab) i \Delta } \left( \CAL{D} \left( 1 + 4ab \right) \phi \right) (z) dz, 
\end{align*}
we have \eqref{K16}.
\end{proof}

\begin{Lem} \label{fac:1}
The following factorization of the propagator holds:
\begin{align*}
e^{i\zeta _1 (s) \Delta /(2 \zeta _2 (s))} e^{i(n-1) \zeta _2 (s) \zeta _2 '(s) |x| ^2 /2} &= (i)^{d/2}e^{i(n-1)\zeta _2 (s) \zeta _2 '(s) |x|^2/(2(1+ (n-1) \zeta _1 (s) \zeta _2 '(s)))} 
 \\ & \quad \times e^{ i \zeta _1 (s) (1+(n-1) \zeta _1(s) \zeta _2 '(s)) \Delta/(2 \zeta _2 (s))} \\ & \qquad \times \CAL{D}\left(1+ (n-1) \zeta _1 (s) \zeta _2 '(s) \right)
\end{align*}
for any $s \geq r_0$.
\end{Lem}
\begin{proof}
The desired identity can be obtained by applying $a = \zeta _1 (s)  /(2 \zeta _2 (s)) $ and $ b=  (n-1) \zeta _2 (s) \zeta _2 '(s)  /2 $ in Lemma \ref{lem:41}. Hence all we need to show is $4ab = (n-1) \zeta_1 (s) \zeta _2 '(s) \neq -1 $. Thanks to the assumption \ref{A1}, as $s \to \infty$, $\zeta _1 (s)$ and $\zeta _2 (s)$ satisfy that  
\begin{align*}
\zeta _1 (s) = c_1 s^{\lambda} + o(s^{\lambda}), \quad \zeta _2 (s) = c_2 s^{1- \lambda} + o(s^{1- \lambda})
\end{align*}
and 
\begin{align*}
\zeta _1 '(s) = c_1 \lambda  s^{\lambda -1} + o(s^{\lambda -1 }), \quad \zeta _2' (s) = c_2(1- \lambda) s^{- \lambda} + o(s^{- \lambda})
\end{align*}
with some non zero constants $c_1$ and $c_2$ (see e.g., Geluk-Mari\'{c}-Tomic \cite{GMT}). Thanks to the condition of discriminant, for all $s$, we have 
\begin{align*}
1 = \zeta _1 (s) \zeta _2'(s) - \zeta _1 '(s) \zeta _2 (s)  = c_1c_2 (1- 2 \lambda) + o(1), 
\end{align*}
which yields 
\begin{align*}
 (n-1) \zeta _1 (s) \zeta _2 '(s) = (n-1) (1-\lambda)c_1c_2 + o(1) = \frac{1-\lambda}{1-2 \lambda} (n-1) + o(1). 
\end{align*}
as $s \to \infty$. If $(n-1) (1- \lambda) /\rho  = -1$ holds, then $\lambda = n/(n+1)$, and which contradicts that $\lambda \in (0,1/2)$.  

\end{proof}


\section{Estimates for non-resonant terms} \label{sec:nre}

In this section, we shall first prove Proposition \ref{nonres}.
Set $\delta>0$ as in \eqref{K15}.
To this end, let us estimate the term associated with $\CAL{N} (u_p(t))$. 
For simplicity, define 
\begin{align*}
	V_{1} (\varphi) ={}& \left\lVert \varphi \right\rVert_{\infty}^{p_c} \left\{ \norm{\varphi}_{H^{\delta'}} \( 1 +  \norm{\varphi}_{\infty}^{p_c} \)^{\lceil \delta \rceil} + \left\| |\xi|^{- \delta} \varphi \right\|_{2} \right\}, \\
	V_{2} (\varphi) ={}& \left\lVert \varphi \right\rVert_{\infty}^{2p_c} \left\{ \norm{\varphi}_{H^{\delta'}}  \( 1 +  \norm{\varphi}_{\infty}^{p_c} \)^{\lceil \delta \rceil} + \left\| |\xi|^{- \delta} \varphi \right\|_{2} \right\}.
\end{align*}
We notice 
\begin{align*}
	u_p(t) = \CAL{M}_1(t) \CAL{D}(\zeta _2 (t)) \widehat{w}(t) 
	= \CAL{D}(\zeta _2 (t)) E(t) \widehat{w} (t), 
\end{align*}
where  
\begin{align*}
E(t) \phi = e^{i \zeta _2 (t) \zeta _2 '(t) x^2/2} \phi (x).
\end{align*}
Then by the same scheme as in \cite{MM2}, one rewrites 
\begin{align*}
	\CAL{N}(u_p(t)) ={}& \sum_{n \neq 0,1} g_n \left| \CAL{D}(\zeta _2 (t)) E(t) \widehat{w}(t) \right|^{1 + p_c -n} \left( \CAL{D}(\zeta _2 (t) ) E(t) \widehat{w}(t) \right)^n \\ 
	={}& \sum_{n \neq 0,1} g_n \left( \frac{c_{dn}}{ \zeta _2 (t)^{1/(1- \lambda)}} \CAL{D}(\zeta _2 (t)) E(t)^n  \phi_n(t) 
\right),
\end{align*}
where $c_{dn} = {i}^{-d(n-1)/2}$ and $\phi_n (t) = \left| \widehat{w}(t)  \right|^{1 + p_c -n} (\widehat{w}(t))^n$.
We here set $\psi_0 \in C^{\infty}_0(\R^d)$ such that 
\[ 
	\psi_0 (y) = 
	\begin{cases}
	1, & |y| \leq 1, \\ 
	0 , & |y| >2,
	\end{cases}
	\quad 0 \leq \psi_0 \leq 1,
\]
and let $\mathcal{K} \coloneqq \mathcal{K}_{\psi_0}(t, n)$ as in \eqref{kpsi:1}.
Note that $\nabla \psi_0 (0) = 0$ and in particular, there are no upper bound of $\theta$ in Lemma \ref{mol:1.1} because of $\psi_0 \equiv 1$ on $[-1, 1]$.
We then decompose $\mathcal{N}(u_p)$ as $\CAL{N}(u_p(t)) = \CAL{P} + \CAL{Q}$
with the low-frequency part
\begin{align*}
\CAL{P} \coloneqq  \sum_{n \neq 0,1} g_n \left( 
\frac{c_{dn}}{ \zeta_2 (t)^{1/(1- \lambda)}} \CAL{D}(\zeta _2 (t)) E(t)^n \CAL{K}  \phi_n(t) 
\right)
\end{align*}
and the high-frequency part
\begin{align*}
\CAL{Q} \coloneqq  \sum_{n \neq 0,1} g_n \left( 
\frac{c_{dn}}{ \zeta_2 (t)^{1/(1- \lambda)}} \CAL{D}(\zeta _2 (t)) E(t)^n \left( 1 - \CAL{K} \right) \phi_n(t) 
\right).
\end{align*}

\subsection{Estimate for $\CAL{Q}$}
The term associated with $\CAL{Q}$ can be estimated as the same argument in \cite{MM2, MMU}. 
Firstly, noting $\nabla \psi_0 (0) = 0$, we see from Lemma \ref{mol:1.1} that
\begin{align*}
	\left\| \int_t^{\infty} U_0 (t,s) \CAL{Q}\, ds \right\|_{\infty, 2, \lambda}
	\leq{}& C \sup_{t \geq \tau} \J{t}^{-\lambda} \int_{t}^{\infty} \sum_{n \neq 0,1} |g_n| \left\|(\CAL{K}-1) \phi_n (s) \right\|_2 \frac{ds}{s} \\
	\leq{}& C \sup_{t \geq \tau} \J{t}^{-\lambda} \int_{t}^{\infty} \sum_{n \neq 0,1} |g_n|  |n|^{-\delta} s^{-1-\frac{\delta}{2}\rho} \left\| |\nabla|^{\delta} \phi_n (s) \right\|_2\, ds \\
	\leq{}& C \J{\tau}^{-\lambda -\frac{\delta}{2}\rho} \J{g_1 \log \tau}^{\lceil \delta \rceil} V_1(\widehat{u_+}) \sum_{n \neq 0,1} |n|^{\gamma - \delta} |g_n|
\end{align*}
for any $\tau \geq r_0$.
Moreover, arguing as in the above, by employing Lemma \ref{Str:1}, one sees from $\lambda - \delta \rho/2 < 0$ that 
\begin{align*}
	\left\| \int_t^{\infty} U_0 (t,s) \CAL{Q}\, ds \right\|_{q, r, \lambda, \tau} 
	\leq{}& C \int_\tau^{\infty} \left\langle t \right\rangle^{\lambda} \left\| \CAL{Q} \right\|_{2}\, dt \\
	\leq{}& C \J{\tau}^{\lambda -\frac{\delta}{2}\rho} \J{g_1 \log \tau}^{\lceil \delta \rceil} V_1(\widehat{u_+}) \sum_{n \neq 0,1} |n|^{\gamma - \delta} |g_n|
\end{align*}
for any $\tau \geq r_0$. 
Hence conclude the estimate
\begin{align}
\begin{aligned}
		&{}\left\| \int_t^{\infty} U_0 (t,s) \CAL{Q}\, ds \right\|_{\infty, 2, \lambda} + \tau^{-2\lambda} \left\| \int_t^{\infty} U_0 (t,s) \CAL{Q}\, ds \right\|_{q, r, \lambda, \tau} \\
		\leq{}& C \J{\tau}^{-\lambda -\frac{\delta}{2}\rho} \J{g_1 \log \tau}^{\lceil \delta \rceil} V_1(\widehat{u_+}) \sum_{n \neq 0,1} |n|^{\gamma - \delta} |g_n|.
\end{aligned}
	\label{est:q}
\end{align}

\subsection{Estimate for $\CAL{P}$}
We estimate the term associated with $\CAL{P}$.  Thanks to Lemma \ref{mdfm:1}, one has 
\begin{align*}
U_0(t,0) = \CAL{M}_1(t) \CAL{D}(\zeta _2 (t)) \SCR{F} \CAL{M} \left( \frac{\zeta _2 (t)}{\zeta _1 (t)} \right) = 
\CAL{M}_1(t) \CAL{D}(\zeta _2 (t))  U \left(- \frac{\zeta _1 (t)}{\zeta _2 (t)} \right) \SCR{F},
\end{align*}
where $U(t) = e^{it \Delta/2}$.
Then $U_0 (t,s) \CAL{P}$ can be rewritten as 
\begin{align*}
U_0(t,s) \CAL{P} &= U_0 (t,0) \SCR{F}^{-1} \SCR{F} U_0 (0,s) \CAL{P} \\ &= 
\CAL{M}_1(t) \CAL{D}(\zeta _2 (t)) U \left( -\frac{\zeta _1 (t)}{ \zeta _2 (t)} \right) \SCR{F} U_0 (0,s).
\end{align*}
We calculate the term $\SCR{F} U_0(0,s) \CAL{P}$. Since 
\begin{align*}
\CAL{D}(\zeta _2 (s)) ^{-1} \CAL{M}_1 (s)^{-1} \CAL{D}(\zeta _2 (s)) \phi &= \frac{1}{(i \zeta _2 (s))^{n/2}}\CAL{D}(\zeta _2 (s)) ^{-1} \left( e^{-i \zeta _2 '(s) x^2/(2 \zeta _2 (s))} \phi (x/\zeta _2 (s)) \right) \\ &= 
e^{-i \zeta _2 '(s) \zeta _2 (s) x^2 /2} \phi(x) = E(s)^{-1} \phi, 
\end{align*}
we have 
\begin{align*}
\SCR{F} U_0(0,s) \CAL{P} ={}& U \left( \frac{\zeta _1 (s)}{\zeta _2 (s)} \right) \CAL{D}(\zeta _2 (s)) ^{-1}\CAL{M}_1(s)^{-1}  \sum_{n \neq 0,1} g_n \left( 
\frac{c_{dn}}{ \zeta _2 (s)^{1/(1- \lambda)}} \CAL{D}(\zeta _2 (s)) E(s)^n \CAL{K}  \phi_n(s) 
\right) \\ 
	={}& U \left( \frac{\zeta _1 (s)}{\zeta _2 (s)} \right) \sum_{n \neq 0,1} g_n \left( 
\frac{c_{dn}}{ \zeta _2 (s)^{1/(1- \lambda)}}  E(s)^{n-1} \CAL{K}  \phi_n(s) \right),
\end{align*}
where $\phi_n (t) = \left| \widehat{w}(t)  \right|^{1 + p_c -n} (\widehat{w}(t))^n$.
Using the lemma \ref{fac:1}, one has 
\begin{align*}
U \left( \frac{\zeta_1 (s)}{\zeta_2 (s)} \right)  E(s)^{n-1} 
	={}& i^{d/2}e^{i(n-1)\zeta_2 (s) \zeta_2^{\prime}(s) |x|^2/(2(1+ (n-1) \zeta_1 (s) \zeta_2 '(s)))} \\ 
	&{} \times e^{ i \zeta _1 (s) (1+(n-1) \zeta_1(s) \zeta_2 '(s)) \Delta/(2 \zeta_2 (s))} 
	\CAL{D}\left(1+ (n-1) \zeta_1 (s) \zeta_2 '(s) \right). 
\end{align*}
Here we set 
\begin{align*}
A_n(s) ={}& \frac{(n-1) \zeta_2 (s) \zeta_2 '(s)}{2(1 + (n-1) \zeta_1 (s) \zeta_2 '(s))}, 
	\quad B_n (s) = \frac{\zeta_1 (s) (1 + (n-1) \zeta_1 (s) \zeta_2 '(s))}{2 \zeta _2 (s)}, \\
	C_n(s) ={}& 1 + (n-1) \zeta_1 (s) \zeta_2 '(s).
\end{align*}
Then it holds that 
\begin{align*}
U \left( \frac{\zeta_1 (s)}{\zeta_2 (s)} \right) E(s)^{n-1} 
	= i^{d/2}e^{iA_n(s) |x|^2} e^{ iB_n(s) \Delta}\CAL{D}\left(C_n(s) \right),
\end{align*}
which implies
\begin{align}
\SCR{F} U_0(0,s) \CAL{D} \( \zeta_2(s) \)  E(s)^{n} 
	= i^{d/2}e^{iA_n(s) |x|^2} e^{ i B_n(s) \Delta}\CAL{D}\left(C_n(s) \right). \label{fac:2}
\end{align}
By using $\zeta _1 (s) \zeta _2 '(s) - \zeta _1 ' (s) \zeta _2 (s) = 1$ for any $s \in {\bf R}$, we compute
\begin{align*}
\frac{d}{ds} A_n(s) ={}& \frac{(n-1) (((\zeta_2 ')^2+ \zeta _2  \zeta _2 '' )(1 + (n-1) \zeta_1  \zeta_2 ') - (n-1)\zeta_2 \zeta_2 ' (\zeta_1 '\zeta_2 ' + \zeta_1 \zeta_2 ''))}{2(1 + (n-1) \zeta_1 (s) \zeta_2 '(s))^2} \\ 
	={}& \frac{{ (n-1) \{ n( \zeta_2'(s) )^2 + \zeta_2 (s) \zeta_2 ''(s) \}}}{2(1 + (n-1) \zeta_1 (s) \zeta_2 '(s))^2}
\end{align*}
and 
\begin{align*}
	\frac{d}{ds} B_n (s) ={}& \frac{d}{ds} \left( \frac{n \zeta_1 (s)}{2 \zeta_2 (s)} + \frac{(n-1) \zeta_1 (s) \zeta_1 '(s)}{2} \right) \\ 
	={}& -\frac{n}{2 \zeta _2 (s)^2} + \frac{n-1}{2} \left( (\zeta _1 '(s)) ^2 + \zeta _1 (s) \zeta _1 ''(s) \right). 
\end{align*}
Hence
\begin{align}
	A_{n}^{\prime}(s) \sim s^{-2 \lambda}, \quad B_{n}^{\prime}(s) \sim n s^{2 \lambda- 2}. \label{coff:1}
\end{align}
Here, for any $\alpha$, $\beta \in \R$, we write $F(t,n) \sim t^{\alpha} n^{\beta}$ if there exist constants $C_1$, $C_2 >0$ such that $C_1 t^{\alpha} |n|^{\beta} \leq F(t,n) \leq C_2 t^{\alpha} |n|^{\beta}$ for any $t \geq r_0$ and $n \neq 0$.
Let 
\begin{align*}
\SCR{A}_n (s) = \left( 1 + is A_n'(s) |x|^2   \right)^{-1}.
\end{align*} 
One has 
\begin{align*}
	&{} \int_t^{\infty} U\left( \frac{\zeta_1 (s)}{\zeta_2 (s)} \right) E(s)^{n-1} \CAL{K} \phi_n (s) \frac{ds}{\zeta _2 (s)^{1/(1- \lambda)} } \\ 
	={}& \int_t^{\infty} \SCR{A}_n (s) \partial _s \left( se^{i A_n (s) |x| ^2} \right) e^{iB_n (s) \Delta} \CAL{D}(C_n (s))\CAL{K} \phi_n (s) \frac{ds}{\zeta _2 (s)^{1/(1- \lambda)} } \\ 
	={}&  - e^{i A_n (t) |x|^2} \SCR{A}_n (t) e^{iB_n (t) \Delta} \CAL{D}(C_n (t))\CAL{K} \phi_n (t) \frac{t}{\zeta _2 (t)^{1/(1- \lambda)} } \\ 
	&{} - \int_t^{\infty} e^{i A_n (s) |x|^2} \partial_s  \left( \frac{\SCR{A}_n (s) }{\zeta_2 (s)^{1/(1- \lambda)} }\right) s  e^{iB_n (s) \Delta} \CAL{D}(C_n (s))\CAL{K} \phi_n (s) {ds} \\ 
	&{} \quad - \int_t^{\infty} e^{i A_n (s) |x|^2} \SCR{A}_n (s)  e^{iB_n (s) \Delta} \left( \partial _s +i B_n'(s) \Delta  \right)\CAL{D}(C_n (s))\CAL{K} \phi_n (s) \frac{s ds}{\zeta _2 (s)^{1/(1- \lambda)} }. 
\end{align*}
The identity \eqref{fac:2} yields
\begin{align*}
	e^{iA_n(s) |x|^2} = i^{-d/2} \SCR{F} U_0(0,s) \CAL{D} \( \zeta_2(s) \)  E(s)^{n} \CAL{D}\left(C_n(s) \right)^{-1} e^{- i B_n(s) \Delta}.
\end{align*}
Hence we see that
\begin{align}
\begin{aligned}
	i \int_t^{\infty} U_0(t,s) \CAL{P}(s) ds  &= i \CAL{M}_1(t) \CAL{D}(\zeta _2 (t)) U\left(- \frac{\zeta _1 (t)}{\zeta _2 (t)} \right) \\ 
	&{} \quad \times\sum_{n \neq 0,1} g_n c_{dn} \int_t^{\infty} U\left( \frac{\zeta _1 (s)}{\zeta _2 (s)} \right) E(s)^{n-1} \CAL{K} \phi_n (s) \frac{ds}{\zeta _2 (s)^{1/(1- \lambda)} } \\
{}& = -i( I_1 + I_2 +I_3 ),
\end{aligned}
	\label{est:p1}
\end{align}
where 
\begin{align*}
I_1 \coloneqq{}& \sum_{n \neq 0,1} g_n i^{-dn/2} \CAL{D} \( \zeta_2(t) \)  E(t)^{n} \CAL{D}\left(C_n(t) \right)^{-1} e^{- i B_n(t) \Delta} \\
	&{} \times \SCR{A}_n (t)  e^{iB_n (t) \Delta} \CAL{D}(C_n (t))\CAL{K} \phi_n (t) \frac{t}{\zeta _2 (t)^{1/(1- \lambda)} }, \\ 
	I_2 \coloneqq{}& \int_t^{\infty} U_0(t,s) \CAL{D} \( \zeta_2(s) \) \sum_{n \neq 0,1} g_n i^{-dn/2}   E(s)^{n} \CAL{D}\left(C_n(s) \right)^{-1} e^{- i B_n(s) \Delta} \\
	&{} \times s \partial_s  \left( \frac{\SCR{A}_n (s) }{\zeta_2 (s)^{1/(1- \lambda)} }\right) e^{iB_n (s) \Delta} \CAL{D}(C_n (s))\CAL{K} \phi_n (s)\, ds,
\end{align*}
and 
\begin{align*}
I_3 \coloneqq{}& \int_t^{\infty} U_0(t,s) \CAL{D} \( \zeta_2(s) \) \sum_{n \neq 0,1} g_n i^{-dn/2}  E(s)^{n} \CAL{D}\left(C_n(s) \right)^{-1} e^{- i B_n(s) \Delta} \\
	&{} \times \frac{s}{\zeta _2 (s)^{1/(1- \lambda)} }\SCR{A}_n (s)  e^{iB_n (s) \Delta} \left( \partial _s +i B_n'(s) \Delta  \right)\CAL{D}(C_n (s))\CAL{K} \phi_n (s)\, ds.
\end{align*}
In order to estimate $I_j$, we need the following:

\begin{Lem} \label{lem:11}
Set $\SCR{G}_n (t) = \SCR{A}_n (t) e^{iB_n (t) \Delta} \CAL{D}(C_n (t))\CAL{K}$. 
Let $\kappa = 2^{-1}\( \delta - d/2\)$.
Take $r>2$.
Then it holds that
\begin{align}
	\Lebn{\SCR{G}_n (t) \phi_n (t)}{2} \leq{}& C t^{-\frac{\delta}{2}\rho } |n|^{-\delta + \kappa} \( \hSobn{\phi_n (t)}{\delta} + |\psi_0(0)| \Lebn{|\xi|^{-\delta} \phi_n (t)}{2} \), 
	\label{lem:12} \\
	\Lebn{e^{-iB_n (t) \Delta} \SCR{G}_n (t) \phi_n (t)}{r}
	\leq{}& C t^{\frac{\nu - \delta}{2}\rho } |n|^{-\delta + \kappa} \( \hSobn{\phi_n (t)}{\delta} + |\psi_0(0)| \Lebn{|\xi|^{-\delta} \phi_n (t)}{2} \) \label{lem:12a}
\end{align} 
for any $t \geq 1$ and $n \neq 0$, $1$, where $\nu = d \(\frac{1}{2} - \frac{1}{r} \)$. 
\end{Lem}

\begin{proof} 
We here remark that
\begin{align}
	A_n(t) \sim t^{\rho}, \quad B_n(t) \sim t^{-\rho}n, \quad C_n(t) \sim n.
\end{align}
Let us first show \eqref{lem:12}.
The triangle inequality gives us
\begin{align*}
	\Lebn{\SCR{G}_n (t) \phi_n (t)}{2} 
	\leq{}& \Lebn{ \SCR{A}_n (t)  \( e^{iB_n (t)\Delta} - 1\) \CAL{D}(C_n (t))\CAL{K} \phi_n (t)}{2} \\
	&{}+ \Lebn{\SCR{A}_n (t) \CAL{D}(C_n (t)) \( \CAL{K} - \psi(0) \) \phi_n (t)}{2} \\
	&{}+ |\psi_0(0)|\Lebn{\SCR{A}_n (t) \CAL{D}(C_n (t)) \phi_n (t)}{2} \\
	=:{}& \mathrm{I}_n + \mathrm{II}_n + \mathrm{III}_n.
\end{align*}
Combining the H\"older inequality and Sobolev embedding with the fact
\[
	|e^{iB_n (t)|\xi|^2} - 1| \leq C |B_n(t)|^{\theta} |\xi|^{2\theta}
\]
for any $\theta \in [0,1]$, it holds that
\begin{align*}
	\mathrm{I}_n \leq{}& C \Lebn{\SCR{A}_n (t)}{p_1}  \Lebn{ |\nabla|^{\frac{n}{p_1}} |B_n(t)|^{\frac{1}{2}\(\delta - \frac{d}{p_1} \)} |\nabla|^{\delta - \frac{d}{p_1}} \CAL{D}(C_n (t))\CAL{K} \phi_n (t)}{2} \\
	\leq{}& C |t A_{n}'(t)|^{-\frac{d}{2p_1}} |B_n(t)|^{\frac{1}{2}\(\delta - \frac{d}{p_1} \)} |C_n(t)|^{-\delta} \Lebn{ |\nabla|^{\delta} \CAL{K} \phi_n (t)}{2} \\
	\leq{}& C t^{-\frac{\delta}{2}\rho} |n|^{-\delta + \frac{1}{2}\(\delta - \frac{d}{p_1} \)} \hSobn{\phi_n (t)}{\delta} 
\end{align*}
for any $p_1 \geq 2$. Note that we are able to choose $p_1$ such that 
\[
	\frac{1}{2}\(\delta - \frac{d}{p_1} \) \leq \frac{1}{2}\(\delta - \frac{d}{2} \) = \kappa. 
\] 
Also, Lemma \ref{mol:1.1} (ii) tells us that for any $p_2 \geq 2$,
\begin{align*}
	\mathrm{II}_n \leq{}& C\Lebn{\SCR{A}_n (t)}{p_2} \Lebn{ |\nabla|^{\frac{d}{p_2}} \CAL{D}(C_n (t)) \( \CAL{K} - \psi_0(0) \) \phi_n (t)}{2} \\
	\leq{}&  C|t A_{n}'(t)|^{-\frac{d}{2p_2}} |C_n (t)|^{-\frac{d}{p_2}} \Lebn{ |\nabla|^{\frac{d}{p_2}} \( \CAL{K} - \psi_0(0) \) \phi_n (t)}{2} \\
	\leq{}&  Ct^{-\frac{d}{2p_2}\rho } |n|^{-\frac{d}{p_2}} |t|^{-\frac{\theta_2}{2}(1-2\lambda	)} |n|^{-\theta_2} \hSobn{\phi_n (t)}{\frac{d}{p_2} + \theta_2} \\
	\leq{}&  C t^{-\frac{\delta}{2}\rho } |n|^{-\delta} \hSobn{\phi_n (t)}{\delta}
\end{align*}
Here we choose $p_2 \geq 2$ and $\theta_2 \geq 0$ such that $\frac{d}{p_2} + \theta_2 = \delta$.
Finally, let us consider the estimate of $\mathrm{III}_n$. 
Note that $|x|^{\theta} |\SCR{A}_n(t)| \leq C (t A_{n}'(t))^{-\frac{\theta}{2}}$ for any $\theta \in [0,2]$.
Therefore it is obtained that 
\begin{align*}
	\mathrm{III}_n \leq{}& C (t A_{n}'(t))^{-\frac{\delta}{2}} |\psi_0(0)| \Lebn{|\xi|^{-\delta} \CAL{D}(C_n (t)) \phi_n (t)}{2} \\
	\leq{}& C (t A_{n}'(t))^{-\frac{\delta}{2}} |C_n (t)|^{-\delta} |\psi_0(0)| \Lebn{|\xi|^{-\delta} \phi_n (t)}{2} \\
	\leq{}& C t^{-\frac{\delta}{2}\rho } |n|^{-\delta} |\psi_0(0)| \Lebn{|\xi|^{-\delta} \phi_n (t)}{2}.
\end{align*}
Collecting these above, we conclude \eqref{lem:12}. Let us move on to the proof of \eqref{lem:12a}.
By using the Gagliardo-Nirenberg inequality, we estimate 
\begin{align}
	&{}\Lebn{\SCR{G}_n (t) \phi_n (t)}{r} \leq C\Lebn{\SCR{G}_n (t) \phi_n (t)}{2}^{1-\nu} \Lebn{\nabla \SCR{G}_n (t) \phi_n (t)}{2}^{\nu},
	\label{eq24a}
\end{align}
where $\nu = d \(\frac{1}{2} - \frac{1}{r} \)$. In terms of $\Lebn{\nabla \SCR{G}_n (t) \phi_n (t)}{2}$, applying 
\[
	|\nabla \SCR{A}_n(t)| \leq C (tA'(t))^{\frac{1}{2}}|\SCR{A}_n(t)|
\]
for any $t \geq 1$, one has
\begin{align}
	\begin{aligned}
	\Lebn{\nabla \SCR{G}_n (t) \phi_n (t)}{2} \leq{}& C(tA_{n}'(t))^{\frac{1}{2}} \Lebn{\SCR{G}_n (t) \phi_n (t)}{2} \\
	&{}+ \Lebn{\SCR{A}_n (t) \nabla e^{iB_n (t) \Delta} \CAL{D}(C_n (t))\CAL{K} \phi_n (t)}{2}.
	\end{aligned}
	\label{eq24}
\end{align}
As for the second term in the above right-hand side,
we here define a regularizing operator by $\widetilde{\mathcal{K}} = \mathcal{K}_{\widetilde{\psi}_0}$ with $\widetilde{\psi}_0 \in C^{\I}_0(\R^d)$ satisfying $\widetilde{\psi}_0(x) = 1$ ($|x| \leq 2$), $\widetilde{\psi}_0(x) = 0$ ($|x| > 2$) and $0 \leq \widetilde{\psi}_0 \leq 1$. It then follows from $\psi_0 \widetilde{\psi}_0 = \psi_0$ that $\widetilde{\CAL{K}} \CAL{K} = \CAL{K}$.
Noting
\begin{align*}
	\CAL{D} (C_n (t)) \widetilde{\mathcal{K}}  = \CAL{D} (C_n (t)) \widetilde{\psi}_0\(\frac{-i \nabla}{|n| t^{\rho/2}} \) 
	= \widetilde{\psi}_0\(\frac{-i C_n(t) \nabla}{|n| t^{\rho/2}} \) \CAL{D} (C_n (t)), 
\end{align*}
we have 
\begin{align*}
	\left\| \SCR{A}_n (t) \nabla e^{iB_n (t)} \CAL{D} (C_n (t)) \CAL{K} \phi_n \right\|_2 
	\leq{}& \sum_{j=1}^d \left\| \SCR{A}_n (t) \partial_j e^{iB_n (t)} \CAL{D} (C_n (t)) \tilde{\CAL{K}}\CAL{K} \phi_n \right\|_2 \\
	\leq{}& \sum_{j=1}^d \left\| \SCR{A}_n (t) \partial_j \widetilde{\psi}_0\(\frac{-i C_n(t) \nabla}{|n| t^{\rho/2}} \) e^{iB_n (t)} \CAL{D} (C_n (t)) \CAL{K} \phi_n \right\|_2.
\end{align*}
Here denote $\widetilde{\mathcal{K}_{j}} = \widetilde{\psi}_j\(\frac{-i C_n(t) \nabla}{|n| t^{\rho/2}} \)$ with $ \widetilde{\psi}_{j}(x) = x_{j} \widetilde{\psi}_0(x)$, which has the relation
\begin{align*}
	\partial_j \widetilde{\psi}_0\(\frac{-i C_n(t) \nabla}{|n| t^{\rho/2}} \) = \frac{i|n|t^{\frac{\rho}{2}}}{C_n(t)} \widetilde{\mathcal{K}}_{j}.
\end{align*}
We compute 
\begin{align*}
	\left\| \SCR{A}_n (t) \partial_j e^{iB_n (t)} \CAL{D} (C_n (t)) \tilde{\CAL{K}}\CAL{K} \phi_n \right\|_2
	\leq \frac{|n| t^{\rho/2}}{|C_n(t)|} \left\| \SCR{A}_n (t) \widetilde{\mathcal{K}}_{j} e^{iB_n (t)} \CAL{D} (C_n (t)) \CAL{K} \phi_n \right\|_2. 
\end{align*}
A commutator calculation shows 
\begin{align*}
\left[ \SCR{A}_n (t), \widetilde{\mathcal{K}}_{j} \right] ={}& \SCR{A}_n (t) \left[ \widetilde{\mathcal{K}}_{j},  itA_n'(t) |\xi|^2 \right] \SCR{A}_n (t) \\ 
	={}& \frac{tA_n'(t) C_n (t) }{|n| t^{\rho/2}} \SCR{A}_n (t) \left(  \xi \cdot (\nabla \widetilde{\mathcal{K}}_{j}) + (\nabla \widetilde{\mathcal{K}}_{j}) \cdot \xi \right) \SCR{A}_n (t) \\ 
	={}& \frac{2 tA_n '(t) C_n (t)}{|n| t^{\rho/2}}  \SCR{A}_n (t) \xi \cdot (\nabla \widetilde{\mathcal{K}}) \SCR{A}_n (t) \\ 
	&{} +\frac{itA_n '(t) (C_n (t))^2}{|n|^2 t^{\rho}} \SCR{A}_n (t) (\Delta \widetilde{\mathcal{K}}) \SCR{A}_n (t),
\end{align*}
where we remark that $\nabla \widetilde{\mathcal{K}}$ and $\Delta \widetilde{\mathcal{K}}$ are defined by
\[
	\nabla \widetilde{\mathcal{K}} = \F^{-1} \nabla \widetilde{\psi}_0 \(\frac{C_n(t) \xi}{|n| t^{\rho/2}} \) \F,
	\quad 
	\Delta \widetilde{\mathcal{K}} = \F^{-1} \Delta \widetilde{\psi}_0 \(\frac{C_n(t) \xi}{|n| t^{\rho/2}} \) \F.
\]
This tells us that 
\begin{align*}
	\left\| \left[ \SCR{A}_n (t), \widetilde{\mathcal{K}}_{j} \right] (\SCR{A}_n (t))^{-1} \right\|_{\mathcal{L}(L^2)} 
	\leq C \frac{ \sqrt{|tA_n '(t)|}}{{t}^{\rho/2}} \leq C.
\end{align*}
From $\SCR{A}_n (t) \widetilde{\mathcal{K}}_{j} = \widetilde{\mathcal{K}}_{j} \SCR{A}_n (t) + \left[ \SCR{A}_n (t) , \widetilde{\mathcal{K}}_{j} \right]$ and $\norm{\widetilde{\mathcal{K}}_{j}}_{\SCR{L} (L^2)} \leq C$, we deduce that
\begin{align*}
	\left\| \SCR{A}_n (t) \nabla e^{iB_n (t)} \CAL{D} (C_n (t)) \CAL{K} \phi_n \right\|_2 
	\leq{}& Ct^{\rho/2} \left\| \SCR{A}_n (t)  e^{iB_n (t)} \CAL{D} (C_n (t)) \CAL{K} \phi_n \right\|_2.
\end{align*}
Together with \eqref{eq24}, one sees from the above that
\[
	\left\| \SCR{A}_n (t) \nabla e^{iB_n (t)} \CAL{D} (C_n (t)) \CAL{K} \phi_n \right\|_2
	\leq C t^{\rho/2} \left\|\SCR{G}_n \phi_n \right\|_2.
\]
Combining \eqref{lem:12} with \eqref{eq24a}, this implies
\begin{align*}
	\Lebn{\SCR{G}_n (t) \phi_n (t)}{r} 
	\leq{}& C t^{\frac{\nu \rho}{2}} \Lebn{\SCR{G}_n (t) \phi_n (t)}{2} \\
	\leq{}& C t^{\frac{\nu - \delta}{2}\rho } |n|^{-\delta + \kappa} \( \hSobn{\phi_n (t)}{\delta} + |\psi(0)| \Lebn{|\xi|^{-\delta} \phi_n (t)}{2} \).
\end{align*}
This completes the proof.
\end{proof}

Let us go back to the estimate for $\CAL{P}$. Let $a(\lambda, \eta) = a(\lambda) + \eta$ and we note that,
by taking $\gamma$ near $\delta$ for $d=3$, we have $-\delta+\kappa+\gamma \leq a(\lambda, \eta)$ since $\kappa = 2^{-1}\( \delta - d/2\) < a(\lambda, \eta)$.

\medskip

\noindent
\underline{\bf Estimation of $I_1$.} \\ 
We first estimate the term $I_1$. Clearly 
\begin{align*}
\left\| 
I_1
\right\|_{\infty, 2, \lambda} \leq C \sum_{n \neq 0,1} |g_n| \sup_{t \geq \tau} \J{t}^{-\lambda} \left\|  \SCR{G}_n (t) \phi_n (t)\right\|_2
\end{align*}
holds. 
By means of Proposition \ref{pro:a2} and \eqref{lem:12}, we deduce that 
\begin{align*}
	\J{t}^{-\lambda} \left\| \SCR{G}_n (t) \phi_n (t)\right\|_2 \leq{}& C \J{t}^{-\frac{\delta}{2}\rho  - \lambda} |n|^{-\delta + \kappa} \( \hSobn{\phi_n (t)}{\delta} + \Lebn{|\xi|^{-\delta} \phi_n (t)}{2} \) \\
	\leq{}& C \J{t}^{-\frac{\delta}{2}\rho  - \lambda} \J{g_1 \log t}^{\lceil \delta \rceil} |n|^{a(\lambda, \eta)} V_1(\widehat{u_+})
\end{align*}
for any $t \geq \tau$, which implies 
\begin{align*}
	\left\| I_1\right\|_{\infty, 2, \lambda} 
	\leq C \J{\tau}^{-\lambda - \frac{\delta}{2} \rho} \J{g_1 \log \tau}^{\lceil \delta \rceil} V_1(\widehat{u_+}) \sum_{n \neq 0,1} |n|^{a(\lambda, \eta)}|g_n|.
\end{align*}

On the other hand, simple calculation shows
\begin{align*}
\left\| \CAL{D}(\zeta _2 (t)) E(t)^n \CAL{D}(C_n (t))^{-1} \Phi \right\|_{r} 
&= |\zeta _2 (t)|^{-d(1/2-1/r) } \left\|E(t)^n \CAL{D}(C_n (t))^{-1} \Phi \right\|_{r} \\ 
& \leq C t^{- 2(1 - \lambda)/q} |n|^{2/q}  \left\|\Phi \right\|_{r}
\end{align*}
for any $t \geq \tau$. Note that $2/q = d \( 1/2 - 1/r \)$.
Using Proposition \ref{pro:a2} and \eqref{lem:12a} again, we obtain
\begin{align*}
\left\| I_1 \right\|_{q, r, \lambda}  
\leq{}& C  \sum_{n \neq 0,1} |n|^{\frac{2}{q}}|g_n| \left( \int_{\tau}^{\infty} \J{t}^{- \lambda - 2(1- \lambda)} \left\| e^{-iB_n (t) \Delta} \SCR{G}_n (t) \phi_n (t) \right\|_{r}^{q} dt \right)^{1/q} \\
\leq{}& C V_1(\widehat{u_+}) \sum_{n \neq 0,1} |n|^{\frac{2}{q} -\delta+\kappa+\gamma}|g_n|
\left( \int_{\tau}^{\infty} \J{t}^{- \frac{\delta}{2} \rho q -1 - \lambda} \J{g_1 \log t}^{\lceil \delta \rceil}\, dt \right)^{1/q} \\
\leq{}& C \J{\tau}^{- \frac{\delta}{2}\rho  - \frac{\lambda}{q}} \J{g_1 \log \tau}^{\lceil \delta \rceil} V_1(\widehat{u_+}) \sum_{n \neq 0,1} |n|^{\frac{2}{q} + a(\lambda, \eta)}|g_n| 
\end{align*}
for any $\tau \geq r_0$, since
\begin{align*}
	- \lambda - 2(1- \lambda) + \frac{2/q - \delta}{2} (1-2 \lambda) q 
	= - \frac{\delta}{2}\rho q -1 - \lambda.
\end{align*}
Therefore, together with $2/q \leq 1$, it holds that
\begin{align}
\begin{aligned}
	&{}\left\lVert I_1 \right\rVert_{\infty, 2, \lambda, \tau} + \J{\tau}^{-2\lambda} \left\lVert I_1 \right\rVert_{q, r, \lambda, \tau} \\
	\leq{}& C \J{\tau }^{-\lambda - \frac{\delta}{2} \rho} \J{g_1 \log \tau}^{\lceil \delta \rceil} V_1(\widehat{u_+}) \sum_{n \neq 0,1} |n|^{1+ a(\lambda, \eta)}|g_n|.
\end{aligned}
	\label{est:p2}
\end{align}

\noindent 
\underline{\bf Estimation of $I_2$.} \\ 
We here estimate $I_2$. By the simple calculation, we get 
\begin{align*}
\partial_s \left( \frac{\SCR{A}_n (s)}{\zeta _2 (s)^{1/(1- \lambda)} }\right) 
	={}& \frac{ \left(\SCR{A}_{n} (s)\right)^{\prime} (\zeta _2 (s))^{1/ (1- \lambda)} - (1- \lambda)^{-1} \zeta _2' (s) \zeta_2 (s)^{\lambda /(1- \lambda)}\SCR{A}_n(s) }{(\zeta _2 (s))^{2/(1- \lambda)}} \\
	={}& (\zeta_{2}(s))^{-1/(1- \lambda)}(\SCR{A}_n(s))' - (1- \lambda)^{-1}(\zeta_{2}(s))^{- \left( 1 + 1/(1- \lambda)\right)}\zeta_{2}'(s)\SCR{A}_{n}(s)
\end{align*}
with 
\begin{align*}
\left( \SCR{A}_{n}(s) \right)^{\prime} ={}& \frac{-i A_n '(s) |x| ^2 -is A''_n(s) |x|^2 }{(1 + is A_n'(s) |x|^2)^2} \\
	={}& -s^{-1} \SCR{A}_{n}(s) + s^{-1} \SCR{A}^{2}(s)- \frac{A^{\prime \prime}(s)}{A^{\prime }(s)} \left( \SCR{A}(s) - \SCR{A}^{2}(s) \right).
\end{align*}
Recall that $\zeta_2 \in C^{3}([r_{0}, \infty))$ and $c |s|^{1-\lambda-j}  < |\zeta_2^{(j)} (s)| < C |s| ^{1- \lambda-j} $ 
with constants $0 < c < C$ for $0 \leq j \leq 3$ and any $|s| \geq r_0$.
Thus $A''_{n}(s) \sim s^{-1-2 \lambda}$, where we use 
\[
	\zeta_2^{(3)}(s) = \(-\sigma(s) \zeta_2(s) \)' = -\sigma'(s) \zeta_2(s) - \sigma(s) \zeta_2'(s) \sim s^{-2-\lambda}.
\]
Therefore, since $|\SCR{A}_{n}(s)| \leq 1$, 
there exists a function $\SCR{B}_{n} (s)$ so that $| \SCR{B}_{n} (s)| \leq C |s| ^{-1}$ and 
\begin{align}
s \partial_s  \left( \frac{\SCR{A}_n (s) }{\zeta_2 (s)^{1/(1- \lambda)} }\right) = \SCR{B}_{n} (s) \SCR{A}_n (s) \label{ba:1}
\end{align}
for any $s \geq r_{0}$.
Then by employing Lemma \ref{Str:1} with an admissible pair $(q, r)$, we obtain
\begin{align*}
\left\| I_2 \right\|_{q, r, \lambda} \leq C \left\| \CAL{D}(\zeta_2 (s)) \sum_{n \neq 0,1} g_n \SCR{G}_{1,n}(s) s \partial_s \left( \frac{\SCR{A}_n (s)}{\zeta_2 (s)^{1/(1- \lambda)} }\right) \SCR{G}_{2,n}(s) \phi_n (s)  \right\|_{1,2,-\lambda}
\end{align*}
with $\SCR{G}_{1,n}(s) = E (s)^n \CAL{D}  (C_n (s))^{-1} e^{-i B_n (s) \Delta}$ and $\SCR{G}_{2,n} (s)  = e^{iB_n (s) \Delta} \CAL{D} (C_n (s)) \CAL{K}$. 
One sees from \eqref{ba:1} that 
\begin{align}
	\left\| I_2 \right\|_{q, r, \lambda} 
	\leq C \sum_{n \neq 0,1} |g_n| \int_{\tau}^{\infty} \J{t}^{\lambda} \left\| \SCR{G}_n(t) \phi_n (t) \right\|_2 \frac{dt}{t}. \label{kes:2}
\end{align}
By Proposition \ref{pro:a2} and \eqref{lem:12}, recalling $\delta >2\lambda/(1-2 \lambda)$, it holds that
\begin{align*}
	\left\| I_2 \right\|_{q, r, \lambda} 
	\leq{}& C V_1(\widehat{u_+})  \sum_{n \neq 0,1} |n|^{a(\lambda, \eta)}|g_n| \int_{\tau}^{\infty} \J{t}^{-1 + \lambda - \frac{\delta}{2} (1-2 \lambda)} \J{g_1 \log t}^{\lceil \delta \rceil}\, dt \\
	\leq{}& C \J{\tau}^{\lambda - \frac{\delta}{2} \rho} \J{g_1 \log \tau}^{\lceil \delta \rceil} V_1(\widehat{u_+}) \sum_{n \neq 0,1} |n|^{a(\lambda, \eta)}|g_n|
\end{align*}
for any $\tau \geq r_0$.
On the other hand, we easily find
\begin{align*}
\left\|  I_2 \right\|_{\infty, 2, \lambda} \leq{}& C \sum_{n \neq 0,1} |g_n|  \sup_{t \geq \tau } \J{t}^{ - \lambda } \int_{t}^{\infty}  \left\| \SCR{G}_n (s) \phi_n (s)  \right\|_2 \frac{ds}{s} \\ 
	\leq{}& C \J{\tau }^{-\lambda - \frac{\delta}{2} \rho} \J{g_1 \log \tau}^{\lceil \delta \rceil} V_1(\widehat{u_+}) \sum_{n \neq 0,1} |n|^{a(\lambda, \eta)}|g_n|
\end{align*}
for any $\tau \geq r_0$.
Hence we conclude that 
\begin{align}
\begin{aligned}
	&{}\left\lVert I_2 \right\rVert_{\infty, 2, \lambda, \tau} + \J{\tau}^{-2\lambda} \left\lVert I_2 \right\rVert_{q, r, \lambda, \tau} \\
	\leq{}& C \J{\tau }^{-\lambda - \frac{\delta}{2} \rho} \J{g_1 \log \tau}^{\lceil \delta \rceil} V_1(\widehat{u_+}) \sum_{n \neq 0,1} |n|^{a(\lambda, \eta)}|g_n|.
\end{aligned}
	\label{est:p3}
\end{align}


\noindent 
\underline{\bf Estimation of $I_3$.} \\ 
Similarly to the estimate of $I_2$, we have
\begin{align}
	\begin{aligned}
	&{}\left\lVert I_3 \right\rVert_{\infty, 2, \lambda, \tau} \\
	\leq{}& C \sum_{n \neq 0,1} |g_n| \sup_{s \geq \tau} \J{s}^{-\lambda} \int_s^{\infty} \left\|\SCR{A}_n (t) e^{iB_n (t) \Delta} \left( \partial_t +i B_n'(t) \Delta \right)\CAL{D}(C_n (t))\CAL{K} \phi_n (t) \right\|_2\, dt
	\end{aligned}
	\label{i3:1}
\end{align}
for any $\tau \geq r_0$. 
Further, Lemma \ref{Str:1} tells us that 
\begin{align}
\begin{aligned}
		&{}\left\lVert I_3 \right\rVert_{q, r, \lambda, \tau} \\
	\leq{}& C \sum_{n \neq 0,1} |g_n| \int_{\tau}^{\infty} \J{t}^{\lambda} \left\|\SCR{A}_n (t)  e^{iB_n (t) \Delta} \left( \partial_t +i B_n'(t) \Delta \right)\CAL{D}(C_n (t))\CAL{K} \phi_n (t) \right\|_2\, dt
\end{aligned}
	\label{i3:2}
\end{align}
for any $\tau \geq r_0$.
Let us here estimate
\[
	\left\|\SCR{A}_n (t) e^{iB_n (t) \Delta} \left( \partial_t +i B_n'(t) \Delta \right)\CAL{D}(C_n (t))\CAL{K} \phi_n (t) \right\|_2.
\]
Firstly we note that
\begin{align*}
	&{}\partial _t \left( \CAL{D} (C_n (t)) \Phi(t, \cdot) \right) (x) \\ 
	={}& \frac{1}{(iC_n (t))^{d/2}} \left(  (\partial_t \Phi ) (t, x/C_n (t)) - \frac{C_n'(t)}{C_n (t)} \left( \frac{d}{2} \Phi (t,x/C_n(t)) + \frac{x}{C_n (t)} \cdot \left( \nabla \Phi \right) (t, x/C_n(t)) \right)\right) \\ 
	={}& \CAL{D}(C_n (t)) \left( \partial _t - \frac{C_n'(t)}{C_n(t)} \left( \frac{d}{2} + x \cdot \nabla \right) \right) \Phi (t, \cdot),
\end{align*}
and 
\begin{align*}
	\Delta \left( \CAL{D} (C_n (t)) \Phi (t, \cdot) \right) 
	= C_n (t)^{-2} \CAL{D} (C_n (t)) \left( \Delta \Phi (t, \cdot) \right).
\end{align*}
These imply
\begin{align*}
&{}\left( \partial _t +i B_n'(t) \Delta  \right)\CAL{D}(C_n (t))\CAL{K} \phi_n (t) \\ 
={}& \CAL{D}(C_n (t)) \left( \partial _t + \frac{i B_n'(t) }{C_n (t) ^2}   \Delta
- \frac{C_n'(t)}{C_n(t)} \left( \frac{d}{2} + x \cdot \nabla \right) 
 \right) \CAL{K} \phi_n (t).
\end{align*}
Here, we define the operators $ \mathcal{K}_{1} = \mathcal{K}_{\psi_{1}}$, $\mathcal{K}_{2} = \mathcal{K}_{\psi_{2}}$, 
where
\[
	\psi_{1}(x) = x \cdot \nabla \psi(x), \quad 
	\psi_{2}(x) = |x|^{2} \psi(x).
\]
Applying the relations
\begin{align}
	\partial _t \CAL{K} = \frac{\rho}{2} t^{-1} \mathcal{K}_{1} + \mathcal{K} \partial_{t}, \quad
	\Delta \mathcal{K} = -|n|^{2} t^{\rho} \mathcal{K}_{2},  
	\label{reop:1}
\end{align}
one obtains
\begin{align*}
	\left( \partial_t +i B_n'(t) \Delta \right)\CAL{D}(C_m (t))\CAL{K} \phi_n (t) 
	={}& \frac{\rho}{2} t^{-1} \CAL{D} (C_n(t)) \mathcal{K}_{1} \phi_{n}(t) + \CAL{D} (C_n(t)) \mathcal{K} \partial_{t} \phi_{n}(t) \\
	&{}-|n|^{2} t^{\rho} \frac{i B_n'(t)}{C_n(t)^2} \mathcal{D}(C_{n}) \mathcal{K}_{2} \phi_{n}(t) \\
	&{}- \frac{d C_n'(t)}{2 C_n(t)} \mathcal{D}(C_{n}) \mathcal{K} \phi_{n}(t) 
	+ \frac{C_n'(t)}{C_n(t)} \mathcal{D}(C_{n}) x \cdot \nabla \mathcal{K} \phi_{n}(t).
\end{align*}
In view of 
\[
	B'_{n}(t) \sim n t^{-1 - \rho}, \quad C_{n}(t) \sim n, \quad C'_{n}(t) \sim n t^{-1},
\]
we then estimate 
\begin{align}
	\begin{aligned}
	&{}\left\| \SCR{A}_m (t)  e^{iB_m (t) \Delta} \left( \partial _t +i B_m'(t) \Delta  \right)\CAL{D}(C_m (t))\CAL{K} \phi_m (t) \right\|_2 \\
	\leq{}& C t^{-1} \norm{\SCR{A}_n (t)  e^{iB_n (t) \Delta} \CAL{D} (C_n(t)) \mathcal{K}_{1} \phi_{n}(t)}_{2} \\
	&{}+ C \norm{\SCR{A}_n (t) e^{iB_n (t) \Delta} \CAL{D} (C_n(t)) \mathcal{K} \partial_{t} \phi_{n}(t)}_{2} \\
	&{}+ C |n| t^{-1} \norm{\SCR{A}_n (t)  e^{iB_n (t) \Delta} \CAL{D} (C_n(t)) \mathcal{K}_{2} \phi_{n}(t)}_{2} \\
	&{}+ C t^{-1} \norm{\SCR{A}_n (t) e^{iB_n (t) \Delta} \CAL{D} (C_n(t)) \mathcal{K} \phi_{n}(t)}_{2} \\
	&{}+ C t^{-1} \norm{\SCR{A}_n (t) e^{iB_n (t) \Delta} \CAL{D} (C_n(t)) x \cdot \nabla \mathcal{K} \phi_{n}(t)}_{2} \\
	\eqqcolon{}& J_1 + J_2  + J_3  + J_4 + J_{5}.
	\end{aligned}
	\label{est:non1}
\end{align}
In terms of $J_{j}$ for $j=1$, $2$, $3$, $4$, since $\mathcal{K}_{1}$ and $\mathcal{K}_{2}$ are of the form \eqref{kpsi:1} and $\psi_{1}(0) = \psi_{2}(0) = 0$, one sees from \eqref{lem:12} and that
\begin{align}
	\begin{aligned}
	J_{1} \leq{}& C t^{-1-\frac{\delta}{2}\rho } |n|^{-\delta + \kappa} \hSobn{\phi_n (t)}{\delta}, \\
	J_{2} \leq{}& C t^{-\frac{\delta}{2}\rho } |n|^{-\delta + \kappa} \( \hSobn{\partial_{t}\phi_n (t)}{\delta} + |\psi(0)| \Lebn{|\xi|^{-\delta} \partial_{t} \phi_n (t)}{2} \), \\
	J_{3} \leq{}& C t^{-1-\frac{\delta}{2}\rho } |n|^{-\delta + 1+ \kappa} \hSobn{\phi_n (t)}{\delta}, \\
	J_{4} \leq{}& C t^{-1-\frac{\delta}{2}\rho } |n|^{-\delta + \kappa} \( \hSobn{\phi_n (t)}{\delta} + |\psi(0)| \Lebn{|\xi|^{-\delta} \phi_n (t)}{2} \).
	\end{aligned}
	\label{est:non2}
\end{align}
Let us deal with $J_{5}$. To get the necessary decay of $t$, we need to take it into account for the regularizing operator to degenerate in the low-frequency region.

\begin{Lem} \label{lemf:1}
Let $\kappa = 2^{-1}\( \delta - d/2\)$.
Then the estimate
\begin{align*}
	\norm{\SCR{A}_n (t) e^{iB_n (t) \Delta} \CAL{D} (C_n(t)) x \cdot \nabla \mathcal{K} \phi(t)}_{2} \leq  C t^{-\frac{\delta}{2}\rho } |n|^{-\delta + 1+ \kappa} \hSobn{\phi (t)}{\delta}.
\end{align*}
holds for any $t \geq 1$. 
\end{Lem}
\begin{proof}[{\bf Proof of Lemma \ref{lemf:1}}]
By using the relation
\begin{align*}
e^{iB_n (t) \Delta} x = \left( x + 2i B_n (t) \nabla \right) e^{iB_n(t) \Delta},
\end{align*}
we see from \eqref{reop:1} that
\begin{align*}
	&{}\norm{\SCR{A}_n (t) e^{iB_n (t) \Delta} \CAL{D} (C_n(t)) x \cdot \nabla \mathcal{K} \phi(t)}_{2} \\
	\leq{}& |C_{n}(t)|^{-1} \norm{\SCR{A}_n (t) e^{iB_n (t) \Delta} x \cdot \CAL{D}(C_n(t)) \nabla \mathcal{K} \phi(t)}_{2} \\
	\leq{}& |C_{n}(t)|^{-1} \norm{\SCR{A}_n (t) \left( x + 2i B_n (t) \nabla \right) e^{iB_m(t) \Delta} \cdot \CAL{D}(C_n(t)) \nabla \mathcal{K} \phi(t)}_{2} \\
	\leq{}& |C_{n}(t)|^{-1} \norm{\SCR{A}_n (t) x e^{iB_n(t) \Delta} \cdot \CAL{D}(C_n(t)) \nabla \mathcal{K} \phi(t)}_{2} \\
	&{}+ 2 |C_{n}(t)|^{-2} |B_n (t)| \norm{\SCR{A}_n (t) e^{iB_n(t) \Delta} \CAL{D}(C_n(t)) \Delta \mathcal{K} \phi(t)}_{2} \\
		\leq{}& C|n|^{-1} \norm{\SCR{A}_n (t) x e^{iB_n(t) \Delta} \cdot \CAL{D}(C_n(t)) \nabla \mathcal{K} \phi(t)}_{2} \\
	&{}+ C |n| |t|^{-1} \norm{\SCR{A}_n (t) e^{iB_n(t) \Delta} \CAL{D}(C_n(t)) \mathcal{K}_{2} \phi(t)}_{2} \\
	\eqqcolon{}& J_{51} + J_{52}.
\end{align*}
Since $J_{52}$ is the same as $J_{3}$, one has
\begin{align*}
	J_{52} \leq C t^{-1-\frac{\delta}{2}\rho } |n|^{-\delta + 1 + \kappa} \hSobn{\phi (t)}{\delta}.
\end{align*}
As for $J_{51}$, it holds that
\begin{align*}
	J_{51} 
	\leq{}& \frac{C|n|^{-1}}{\sqrt{|tA'_m(t) |}} \sum_{j=1}^d \left\| |\SCR{A}_n (t)|^{1/2} e^{iB_n(t) \Delta} \CAL{D} (C_n (t)) \partial_j \CAL{K} \phi(t) \right\|_2.
\end{align*}
We here define the regularized operator by $\mathcal{K}_{3j} = \mathcal{K}_{\psi_{3j}}$ with $ \psi_{3j}(x) = x_{j} \psi_0(x)$, which has the relation
\begin{align}
	\partial_{j} \mathcal{K} = i|n|t^{\frac{\rho}{2}} \mathcal{K}_{3j}. \label{reop:2}
\end{align}
Noting $\psi_{3j}(0) =0$, arguing as in Lemma \ref{lem:11}, one estimates
\begin{align*}
	&{}\left\| |\SCR{A}_n (t)|^{1/2} e^{iB_n(t) \Delta} \CAL{D} (C_n (t)) \CAL{K}_{3j} \phi(t) \right\|_2 \\
	\leq{}& \Lebn{ |\SCR{A}_n (t)|^{1/2} \( e^{iB_n (t)\Delta} - 1\) \CAL{D}(C_n (t))\CAL{K}_{3j} \phi (t)}{2} \\
	&{}+ \Lebn{|\SCR{A}_n (t)|^{1/2} \CAL{D}(C_n (t)) \( \CAL{K}_{3j} - \psi_{3j}(0) \) \phi (t)}{2} \\
	=:{}& \mathrm{IV} + \mathrm{V}.	
\end{align*}
Therefore, it is deduced that
\begin{align*}
	\mathrm{IV} \leq{}& C \Lebn{\SCR{A}_n (t)}{p_3/2}^{1/2} \Lebn{ |\nabla|^{\frac{n}{p_3}} |B_n(t)|^{\frac{1}{2}\(\delta - \frac{d}{p_3} \)} |\nabla|^{\delta - \frac{d}{p_3}} \CAL{D}(C_n (t))\CAL{K}_{3j} \phi (t)}{2} \\
	\leq{}& C t^{-\frac{\delta}{2}\rho} |n|^{-\delta + \frac{1}{2}\(\delta - \frac{d}{p_3} \)} \hSobn{\phi (t)}{\delta} 
\end{align*}
for any $p_3 > d$. We remark that we are able to choose $p_3$ such that 
\[
	\frac{1}{2}\(\delta - \frac{d}{p_3} \) = \frac{1}{2}\(\delta - \frac{d}{2} \) 
	+ \frac{d}{2} \left( \frac{1}{2} - \frac{1}{p_{3}} \right) < \kappa +1.
\]
Also, it follows from the property of $\CAL{K}_{3j}$ that for any $p_4 \geq 2$,
\begin{align*}
	\mathrm{V} \leq{}& C\Lebn{\SCR{A}_n (t)}{p_4/2}^{1/2} \Lebn{ |\nabla|^{\frac{d}{p_4}} \CAL{D}(C_n (t)) \( \CAL{K}_{3j} - \psi_{3j}(0) \) \phi (t)}{2} \\
	\leq{}&  Ct^{-\frac{d}{2p_4}\rho } |n|^{-\frac{d}{p_4}} |t|^{-\frac{\theta_3}{2}\rho} |n|^{-\theta_3} \hSobn{\phi (t)}{\frac{d}{p_4} + \theta_2} \\
	\leq{}&  C t^{-\frac{\delta}{2}\rho } |n|^{-\delta} \hSobn{\phi (t)}{\delta}
\end{align*}
Here we choose $p_4 > d$ and $\theta_3 \in [0,1]$ such that $\frac{d}{p_4} + \theta_3 = \delta$.
Combining these estimates with \eqref{reop:2}, together with $t A'_{n}(t) \sim t^{\rho}$, one has
\begin{align*}
	J_{51} \leq{}& C \sum_{j=1}^{d} \left\| |\SCR{A}_n (t)|^{1/2} e^{iB_n(t) \Delta} \CAL{D} (C_n (t)) \CAL{K}_{3j} \phi(t) \right\|_2 
	\leq C t^{-\frac{\delta}{2}\rho } |n|^{-\delta + 1+ \kappa} \hSobn{\phi (t)}{\delta} 
\end{align*}
Hence we conclude that
\[
	\norm{\SCR{A}_n (t) e^{iB_n (t) \Delta} \CAL{D} (C_n(t)) x \cdot \nabla \mathcal{K} \phi(t)}_{2} 
	\leq C t^{-\frac{\delta}{2}\rho } |n|^{-\delta + 1+ \kappa} \hSobn{\phi (t)}{\delta}.
\]
\end{proof}

Let us go back to estimate $J_{5}$. Lemma \ref{lemf:1} leads to 
\[
	J_{5} \leq C t^{-1-\frac{\delta}{2}\rho } |n|^{-\delta + 1+ \kappa} \hSobn{\phi_n (t)}{\delta}.
\]
Collecting the above, \eqref{est:non1} and \eqref{est:non2}, we see from Proposition \ref{pro:a2} and Proposition \ref{pro:a4} that
\begin{align*}
	&{}\left\| \SCR{A}_n (t)  e^{iB_n (t) \Delta} \left( \partial _t +i B_n'(t) \Delta  \right)\CAL{D}(C_n (t))\CAL{K} \phi_n (t) \right\|_2 \\
	\leq{}& C t^{-\frac{\delta}{2}\rho } |n|^{-\delta + \kappa} \( \hSobn{\partial_{t}\phi_n (t)}{\delta} + \Lebn{|\xi|^{-\delta} \partial_{t} \phi_n (t)}{2} \) \\
	&{}+ C t^{-1-\frac{\delta}{2}\rho } |n|^{-\delta + 1+ \kappa} \( \hSobn{\phi_n (t)}{\delta} + \Lebn{|\xi|^{-\delta} \phi_n (t)}{2} \) \\
	\leq{}& C t^{-1-\frac{\delta}{2}\rho } \J{g_1 \log t}^{\lceil \delta \rceil} |n|^{1 -\delta + \kappa + \gamma} \( V_1 (\widehat{u_+}) + V_2(\widehat{u_+}) \).
\end{align*}
Here note that $-\delta + \kappa + \gamma \leq a(\lambda, \eta)$.
Therefore, from \eqref{i3:1} and \eqref{i3:2}, together with $\lambda < \delta \rho/2$, it is obtained that 
\begin{align}
\begin{aligned}
	&{}\left\lVert I_3 \right\rVert_{\infty, 2, \lambda, \tau} + \J{\tau}^{-2\lambda} \left\lVert I_3 \right\rVert_{q, r, \lambda, \tau} \\
	\leq{}& C \J{\tau }^{-\lambda - \frac{\delta}{2} \rho} \J{g_1 \log \tau}^{\lceil \delta \rceil} \( V_1(\widehat{u_+}) + V_2 (\widehat{u_+}) \) \sum_{n \neq 0,1} |n|^{1+ a(\lambda, \eta)}|g_n|.
\end{aligned}
	\label{est:p4}
\end{align}
In conclusion, Proposition \ref{nonres} is established from \eqref{est:p1}, \eqref{est:p2}, \eqref{est:p3} and \eqref{est:p4}.


\subsection{Proof of Main results}

\begin{proof}[Proof of Theorem \ref{T1}]
Take $b_0 >0$ as in \eqref{Mb:1}.
Combing Lemma \ref{res:1} with Lemma \ref{mca:1} and Proposition \ref{nonres}, Proposition \ref{prop:31} tells us that 
there exists a constant $\varepsilon_0 = \varepsilon_0(\norm{g}_{\mathrm{Lip}})>0$ such that
if $\norm{\widehat{u_+}}_{L^\I} \leq \varepsilon_0$, then there exists a $T_{0} \coloneqq T(b_0) \geq r_0$ such that our integral equation \eqref{inteq:1} admits an unique solution $u \in C([T_0, \infty)\, ;\, L^2)$ with $\norm{u}_{X_{T_0, b_0}} \leq M$ for some $M>0$.
Let us show that the solution $u$ satisfies $\norm{u}_{X_{T_0, b}}< \infty$ for any $b \in \(2 \lambda, \lambda + \delta (1-2\lambda)/2 \)$.
Since $\norm{u}_{X_{T_{0}, b}}< \infty$ for any $b \in (2\lambda, b_0]$ is trivial, we will prove $\norm{u}_{X_{T_{0}, b}}< \infty$ for all $b \in (b_0, \lambda + \delta (1-2\lambda)/2)$. Fix $b_1 \in (b_0, \lambda + \delta (1-2\lambda)/2)$. 
Applying Proposition \ref{prop:31} again, there exists a $T_{1} = T(b_1) \geq r_0$ such that there exists a unique solution $u_1 \in C([T_1, \infty)\, ;\, L^2)$ with $\norm{u_1}_{X_{T_1, b_1}} \leq M$.
Without loss of generality, we assume $T_{1} \geq T_{0}$.
One sees from $b_1 > b_0$ that for any $t \geq T_1 \geq 1$,
\begin{align*}
	t^{b_0-\lambda} \norm{u_1 - u_p}_{\infty, 2, \lambda, t} + t^{b_0-2\lambda} \norm{u_1 - u_p}_{q, r, \lambda, t} 
	\leq {T_1}^{b_0 -b_1} \norm{u_1}_{X_{T_1, b}} \leq M,
\end{align*}
which implies $\norm{u_1}_{X_{T_1, b_0}} \leq M$. Hence $u$, $u_1 \in X_{T_{1}, b_0, M}$.
By the uniqueness property of $X_{T_{1}, b_0,M}$, 
it holds that $u = u_1$ on $ [T_1, \infty)$.
Hence we can extend the existence time of $u_1$ as $T_{1} = T_{0}$, that is, $T_{0}$ does not depend on any $b$. Denoting this $T_0$ by $T$, we then have $\norm{u}_{X_{T, b}} < \infty$ for any $b \in (2\lambda, \lambda + \delta (1-2\lambda)/2)$.
Combining an estimate as in \eqref{pr:31b} with the standard argument of the uniqueness, 
the assertion for the uniqueness stated in Remark \ref{rem:18} is immediate.
Moreover, the estimate \eqref{main:2} follows from Proposition \ref{prop:31}.
Hence, the proof is completed.
\end{proof}


\section*{Appendix}
In this appendix, we shall prove the generalized version of Proposition \ref{pro:a1} - Proposition \ref{pro:a4}, replacing $\widehat{w}$ by the generalized form $\phi \exp(i \eta |\phi|^{\alpha})$.
From now on, we set $\eta$, $\mu \in \R$ and let $\alpha \geq 1$ if $d=1$, $2$ and $\alpha \in (1/2, 1)$ if $d=3$. 
Other notations are same as in the beginning of section \ref{sec:3}.

\begin{Prop} \label{pro:ga1}
Let $d \leq 3$.
There exists $C >0$ such that
\begin{align*}
\norm{\phi \exp(i \eta |\phi|^{\alpha})}_{H^{\delta}} \leq{}& C \norm{\phi}_{H^{\delta}}\( 1 + |\eta| \norm{\phi}_{\infty}^{\alpha} \)^{\lceil \delta \rceil} 
\end{align*}
for any $\phi \in H^{\delta}$. 
\end{Prop}

\begin{Prop} \label{pro:ga2}
Let $d \leq 3$.
Set $\Psi = \phi \exp (i\eta |\phi|^{\alpha})$ with $\phi \in H^{\delta'}$.
Then there exists $C >0$ such that
\begin{align*}
\norm{|\Psi|^{\alpha+1-n} \Psi^n}_{H^{\delta}} \leq{}& C \left\langle n \right\rangle^{\gamma} \J{\eta}^{\lceil \delta \rceil} \left\lVert \phi \right\rVert_{\infty}^{\alpha} \norm{\phi}_{H^{\delta'}}\( 1 + \norm{\phi}_{\infty}^{\alpha} \)^{\lceil \delta \rceil}.
\end{align*}
\end{Prop}
When $\Psi = \phi \exp (i \nu |\phi|^{\alpha} \log t)$, a direct computation shows
\begin{align*}
	\partial_t \left(  |\Psi|^{\alpha+1-n} \Psi^n \right) = \frac{i\nu}{t} n |\Psi|^{2\alpha+1-n} \Psi^n.
\end{align*}
Noting $2\alpha \geq 1$, arguing as in Proposition \ref{pro:ga2}, we easily verify the following:

\begin{Prop} \label{pro:ga4}
Set $\Psi = \phi \exp (i \nu |\phi|^{\alpha} \log t)$ with $\phi \in H^{\delta,0}$.
Then there exists $C >0$ such that
\begin{align*}
	\norm{\partial_t \left( |\Psi|^{\alpha+1-n} \Psi^n  \right)}_{H^{\delta}} \leq{}& C |\nu| t^{-1} \J{n}^{1+\delta} \J{\nu \log t}^{\lceil \delta \rceil} \norm{\phi}_{\infty}^{2\alpha} \norm{\phi}_{H^{\delta'}}\( 1 + \norm{\phi}_{\infty}^{\alpha} \)^{\lceil \delta \rceil}.
\end{align*}
\end{Prop}

In order to prove the above Propositions, we need the following lemmas:
\begin{Lem} \label{lem:nl1}
Let $\alpha>0$ and $n \in \Z$. If $\alpha \geq 1$, then
\begin{align}
	\left| |z|^{\alpha-n}z^n - |w|^{\alpha-n}w^n \right| \leq C |n| \left(|z|^{\alpha-1} + |w|^{\alpha-1}\right) |z-w| \label{nl:1}
\end{align}
for any $z$, $w \in \C$. Further, when $ \alpha \in (0,1)$, the following holds:
\begin{align}
	\left| |z|^{\alpha-n}z^n - |w|^{\alpha-n}w^n \right| \leq C |n|^{\alpha} |z-w|^{\alpha}. \label{nl:2}
\end{align}
\end{Lem}
\begin{proof}
\eqref{nl:1} is immediate from the standard calculation.
\eqref{nl:2} is given by \cite[Lemma 2.4]{GV1}.
\end{proof}

\begin{Lem}[\cite{KaP, GraO}]
\label{lem_lei}
Let $s>0$, $1<r<\infty$, $1<p_1,p_2,q_1,q_2\leq \infty$ and $1/r=1/p_1+q_1=1/p_2+1/q_2$. Then we have the following fractional Leibniz rule: 
$$
\norm{|\nabla|^s(fg)}_{L^r} \leq C \( \norm{|\nabla|^sf}_{L^{p_1}}\norm{g}_{L^{q_1}}+\norm{f}_{L^{p_2}}\norm{|\nabla|^sg}_{L^{q_2}}\).
$$
\end{Lem}

\begin{Lem}[\cite{Vi1}]
\label{lem_cha}
Let $F$ be a H\"older continuous function of order  $0< \rho <1$. Suppose that $0<\sigma<\rho$, $1<p<\infty$ and $\sigma/\rho<s<1$. Then
\begin{align*}
\norm{|\nabla|^\sigma F(f)}_{L^p} \leq C \norm{f}_{L^{\infty}}^{\rho-\sigma/s}\norm{|\nabla|^sf}_{L^{p\sigma /s}}^{\sigma/s}.
\end{align*}
\end{Lem}

\begin{Lem} \label{lem:g3a}
When $d = 2$, $3$, the estimate
\begin{align*}
	\norm{|\nabla|^{\delta-1} \left( |\phi|^{\alpha-n}\phi^n \right)}_{\frac{2\delta}{\delta-1}} \leq{}& C \left\langle n \right\rangle^{\gamma-1} \norm{\phi}^{\alpha + 1-\delta}_{\infty} \norm{\nabla \phi}_{2\delta}^{\delta-1}.
\end{align*}
holds. Moreover, if $d=2$, then 
\begin{align*}
	\norm{|\nabla|^{\delta-1} \exp \(i \eta |\phi|^{\alpha}\)}_{\frac{2\delta}{\delta-1}} 
	\leq{}& C |\eta|^{\delta-1} \norm{\phi}^{(\alpha -1) (\delta-1)}_{\infty} \norm{\nabla \phi}_{2\delta}^{\delta-1}.
\end{align*}
When $d=3$, the following holds:
\begin{align*}
	\norm{|\nabla|^{\delta-1} \exp \(i \eta |\phi|^{\alpha}\)}_{\frac{2\delta}{\delta-1}} 
	\leq{}& C |\eta|^{\frac{\delta-1}{\gamma-1}} \norm{\phi}^{\left( \frac{\alpha}{\gamma-1} -1 \right)  (\delta-1)}_{\infty} \norm{\nabla \phi}_{2\delta}^{\delta-1}.
\end{align*}
\end{Lem}
\begin{proof}[Proof of Lemma \ref{lem:g3a}]
Firstly, we shall show the first estimate in $d=2$.  
By $|\nabla \left( |\phi|^{\alpha-n}\phi^n \right)| \leq C \left\langle n \right\rangle |\phi|^{\alpha-1} |\nabla \phi|$, the Gagliardo-Nirenberg inequality implies
\begin{align}
	\begin{aligned}
	\norm{|\nabla|^{\delta-1} \left( |\phi|^{\alpha-n}\phi^n \right)}_{L^{\frac{2\delta}{\delta-1}}} 
	\leq{}& C \norm{|\phi|^{\alpha-n}\phi^n}_{\infty}^{2-\delta} \norm{\nabla \left( |\phi|^{\alpha-n}\phi^n \right)}_{2\delta}^{\delta-1} \\
	\leq{}& C \left\langle n \right\rangle^{\delta-1} \norm{\phi}^{\alpha + 1 -\delta}_{\infty} \norm{\nabla \phi}_{2\delta}^{\delta-1},
	\end{aligned}
	\label{app:n1}
\end{align}
where we note the identity 
\[
	\frac{\delta-1}{2\delta} - \frac{\delta-1}{d} = \( \frac{1}{2\delta} - \frac{1}{d} \) (\delta-1).
\]
This completes the proof of the first estimate in $d= 2$. 
Let us prove the case $d=3$. 
The Gagliardo-Nirenberg inequality gives us 
\begin{align}
	\norm{|\nabla|^{\delta-1} \left( |\phi|^{\alpha-n}\phi^n \right)}_{\frac{2\delta}{\delta-1}} 
	\leq{}& \left\lVert |\phi|^{\alpha-n}\phi^n \right\rVert_{\infty}^{1-(\delta-1)/s_0} \left\lVert |\nabla|^{s_0} \left( |\phi|^{\alpha-n}\phi^n \right) \right\rVert_{2\delta/s_0}^{(\delta-1)/s_0},
	\label{app:n2}
\end{align}
where $s_0 = \alpha - \varepsilon$ and $\varepsilon \in (0, \alpha-\delta+1)$ is chosen later.
Here $|\phi|^{\alpha-n}\phi^n$ are $\alpha$-H\"older continuous functions, since by \eqref{nl:2} it holds that
\begin{align*}
	||z_1|^{\alpha-n}z_1^n - |z_2|^{\alpha-n}z_2^n| \leq C |n|^{\alpha} |z_1-z_2|^{\alpha}.
\end{align*}
Hence, a use of Lemma \ref{lem_cha} yields 
\begin{align*}
	\norm{|\nabla|^{s_0} \left( |\phi|^{\alpha-n}\phi^n \right)}_{2\delta/s_0} 
	\leq{}& C |n|^{\alpha} \norm{\phi}_{\infty}^{\alpha - \frac{s_0}{s}} \norm{|\nabla|^s \phi}_{\frac{2\delta}{s}}^{\frac{s_0}{s}}
\end{align*}
for any $s \in \(s_0/\alpha, 1\)$.
Applying the Gagliaro-Nirenberg inequality, we have 
\[
	\norm{|\nabla|^s \phi}_{\frac{2\delta}{s}} \leq C \norm{\phi}_{\infty}^{1-s} \norm{\nabla \phi}_{2\delta}^s,
\]
which implies
\begin{align*}
	\norm{|\nabla|^{s_0} \left( |\phi|^{\alpha-n}\phi^n \right)}_{2\delta/s_0} 
	\leq{}& C  |n|^{\alpha} \norm{\phi}_{\infty}^{\alpha-s_0} \norm{\nabla \phi}_{2\delta}^{s_0}.
\end{align*}
Thus we see that 
\begin{align*}
	\norm{|\nabla|^{\delta-1} \left( |\phi|^{\alpha-n}\phi^n \right)}_{\frac{2\delta}{\delta-1}} 
	\leq{}& C |n|^{\frac{\alpha}{\alpha-\varepsilon}(\delta-1)} \norm{\phi}^{\alpha + 1-\delta}_{\infty} \norm{\nabla \phi}_{2\delta}^{\delta-1}.
\end{align*}
Taking $\varepsilon>0$ small such that $\frac{\alpha}{\alpha-\varepsilon} (\delta-1) \leq \gamma-1$, the desired estimate is established.

Let us show the second estimate.
Set $\psi = \exp \(i \eta |\phi|^{\alpha}\)$.
In light of $|\nabla \psi| \leq C |\eta| |\phi|^{\alpha-1} |\nabla \phi|$, arguing as in \eqref{app:n1}, we have
\begin{align*}
	\norm{|\nabla|^{\delta-1} \psi}_{\frac{2\delta}{\delta-1}} 
	\leq{}& C \norm{\psi}_{\infty}^{2-\delta} \norm{\nabla \psi}_{2\delta}^{\delta-1} 
	\leq C |\eta|^{\delta-1} \norm{\phi}^{(\alpha -1)(\delta-1)}_{\infty} \norm{\nabla \phi}_{2\delta}^{\delta-1}
\end{align*}
as desired in $d=2$. 
When $d=3$, in the same manner of \eqref{app:n2}, one has
\begin{align*}
	\norm{|\nabla|^{\delta-1} \psi}_{\frac{2\delta}{\delta-1}} 
	\leq C \left\lVert |\nabla|^{s_0} \psi \right\rVert_{2\delta/s_0}^{(\delta-1)/s_0}.
\end{align*}
Noting $|\psi(z_1) - \psi(z_2)| \leq C |\eta| |z_1-z_2|^{\alpha}$, where $\psi(z) = \exp \(i \eta |z|^{\alpha}\)$, it follows from Lemma \ref{lem_cha} and the Gagliaro-Nirenberg inequality that
\begin{align*}
	\norm{|\nabla|^{s_0} \psi}_{2\delta/s_0} 
	\leq{}& C |\eta| \norm{\phi}_{\infty}^{\alpha - \frac{s_0}{s}} \norm{|\nabla|^s \phi}_{\frac{2\delta}{s}}^{\frac{s_0}{s}} \\
	\leq{}& C |\eta| \norm{\phi}_{\infty}^{\alpha - \frac{s_0}{s}} \left( \norm{\phi}_{\infty}^{1-s} \norm{\nabla \phi}_{2\delta}^s  \right)^{\frac{s_0}{s}} \\
	\leq{}& C |\eta| \norm{\phi}_{\infty}^{\alpha-s_0} \norm{\nabla \phi}_{2\delta}^{s_0}.
\end{align*}
for any $s \in \(s_0/\alpha, 1\)$. 
Hence, we conclude 
\begin{align*}
	\norm{|\nabla|^{\delta-1} \psi}_{\frac{2\delta}{\delta-1}} 
	\leq{}& C |\eta|^{\frac{\delta-1}{\gamma-1}} \norm{\phi}^{\left( \frac{\alpha}{s_0} -1 \right)  (\delta-1)}_{\infty} \norm{\nabla \phi}_{2\delta}^{\delta-1},
\end{align*}
taking $\gamma-1 = s_0$.
This completes the proof.
\end{proof}

\begin{proof}[Proof of Proposition \ref{pro:ga1}]
Set
\[
	\psi (\phi)= \exp \(i \eta |\phi|^{\alpha}\).
\]
Let us first show the case $d=1$. Lemma \ref{lem_lei} gives us 
\begin{align*}
	\norm{\phi \exp(i \eta |\phi|^{\alpha})}_{\dot{H}^{\delta}} \leq{}& C \( \norm{\phi}_{\dot{H}^{\delta}} \norm{\psi}_{L^{\infty}} + \norm{\phi}_{L^{\infty}} \norm{\psi}_{\dot{H}^{\delta}} \) \\
	\leq{}& C \( \norm{\phi}_{\dot{H}^{\delta}} + \norm{\phi}_{L^{\infty}} \norm{\psi}_{\dot{H}^{\delta}} \).
\end{align*}
Using \eqref{nl:1}, we see that for any $z_1$, $z_2 \in \mathbb{R}$,
\begin{align*}
	|\psi(z_1) - \psi(z_2)| 
	\leq{}& C|\eta| \left( |z_1|^{\alpha-1} + |z_2|^{\alpha-1} \right)|z_1 - z_2|,
\end{align*}
since $\alpha \geq 1$. By the representation of $\dot{H}^{\delta}$-norm via the Gagliardo semi-norm (see \cite[Proposition 3.4]{DPV}), one sees that
\begin{align*}
	\norm{\psi}_{\dot{H}^{\delta}} \leq{}& C \(\iint_{\R^n \times \R^n} \frac{|\psi(x) - \psi(y)|^2}{|x-y|^{d+2\gamma}}\, dxdy \)^{1/2} 
	\leq C |\eta| \norm{\phi}_{L^{\infty}}^{\alpha-1} \norm{\phi}_{\dot{H}^{\delta}}.
\end{align*}
Combining these above, it holds that
\begin{align*}
\norm{\phi \exp(i \eta |\phi|^{\alpha})}_{\dot{H}^{\delta}} \leq{}& C \norm{\phi}_{\dot{H}^{\delta}} 
	\( 1 + |\eta| \norm{\phi}_{L^{\infty}}^{\alpha} \).
\end{align*}
Since the estimate $\norm{\phi \exp(i \eta |\phi|^{\alpha})}_{L^2} \leq \norm{\phi}_{L^2}$, we have desired estimate.

We shall deal with the case $d= 2$, $3$.
One computes
\begin{align*}
	\nabla \(\phi \psi \) 
	={}& \nabla \phi \psi + i \eta \alpha \( |\phi|^{\alpha-2} \phi^2 \overline{\nabla \phi} + |\phi|^{\alpha} \nabla \phi\) \psi. 
\end{align*}
Denote $F_1 (\phi) = |\phi|^{\alpha-2} \phi^2$ and $F_2(\phi) = |\phi|^{\alpha}$.
Applying Lemma \ref{lem_lei}, we see from the H\"older inequality that
\begin{align*}
	\norm{\phi \exp(i \eta |\phi|^{\alpha})}_{\dot{H}^{\delta}} 
	\leq{}& C \norm{\phi}_{\dot{H}^{\delta}} + \norm{\nabla \phi}_{2\delta} \norm{|\nabla|^{\delta-1} \psi}_{\frac{2\delta}{\delta-1}} \\
	&{}+ C |\eta| \( \norm{|\nabla|^{\delta-1} F_1(\phi)}_{\frac{2\delta}{\delta-1}} + \norm{|\nabla|^{\delta-1} F_2(\phi)}_{\frac{2\delta}{\delta-1}} \) \norm{\nabla \phi}_{2\delta} \\
	&{}+ C|\eta| \norm{\phi}_{L^{\infty}}^{\alpha} \norm{\phi}_{\dot{H}^{\delta}} 
	+ C|\eta| \norm{\phi}_{L^{\infty}}^{\alpha} \norm{\nabla \phi}_{2\delta} \norm{|\nabla|^{\delta-1}\psi}_{\frac{2\delta}{\delta-1}}.
\end{align*}
When $d=2$, combing Lemma \ref{lem:g3a} with the Gagliardo-Nirenberg inequality 
\begin{align}
	\norm{\nabla \phi}_{2\delta} \leq C \norm{\phi}_{\infty}^{1-\frac{1}{\delta}} \norm{\phi}_{\dot{H}^{\delta}}^{\frac{1}{\delta}},
	\label{gn:1}
\end{align}
one obtains
\begin{align*}
	\norm{\phi \exp(i \eta |\phi|^{\alpha})}_{H^{\delta}} 
	\leq{}& C \norm{\phi}_{H^{\delta}} + |\eta|^{\delta-1} \norm{\phi}^{\alpha(\delta-1)}_{\infty} \norm{\phi}_{\dot{H}^{\delta}} \\
	&{}+ C|\eta| \norm{\phi}_{\infty}^{\alpha} \norm{\phi}_{\dot{H}^{\delta}} + C|\eta|^{\delta} \norm{\phi}_{\infty}^{\alpha \delta} \norm{\phi}_{\dot{H}^{\delta}} \\
	\leq{}& C \norm{\phi}_{H^{\delta}} \(1 + |\eta| \norm{\phi}_{\infty}^{\alpha} \)^{\delta}.
\end{align*}
When $d=3$, combing Lemma \ref{lem:g3a} with \eqref{gn:1}, we conclude
\begin{align*}
	\norm{\phi \exp(i \eta |\phi|^{\alpha})}_{H^{\delta}} 
	\leq{}& C \norm{\phi}_{H^{\delta}} \( 1 + |\eta|^{\frac{\delta-1}{\gamma-1}} \norm{\phi}^{\frac{\alpha}{\gamma-1}(\delta-1)}_{\infty} + |\eta| \norm{\phi}_{\infty}^{\alpha} + |\eta|^{1+\frac{\delta-1}{\gamma-1}} \norm{\phi}_{\infty}^{\left({1+\frac{\delta-1}{\gamma-1}} \right)\alpha} \) \\
	\leq{}& C \norm{\phi}_{H^{\delta}} \(1 + |\eta| \norm{\phi}_{\infty}^{\alpha} \)^{1+\frac{\delta-1}{\gamma-1}}.
\end{align*}
This completes the proof.
\end{proof}

\begin{proof}[Proof of Proposition \ref{pro:ga2}]
Let us first show the case $d=1$.
We see from the Gagliardo-Nirenberg inequality that
\begin{align*}
	\norm{|\Psi|^{\alpha+1-n} \Psi^n}_{\dot{H}^{\delta}} 
	\leq{}& C \left\lVert \phi \right\rVert_{\infty}^{\alpha(1-\delta)} \left\lVert \phi \right\rVert_{2}^{1-\delta} \left\lVert |\Psi|^{\alpha+1-n} \Psi^n \right\rVert_{\dot{H}^{1}}^{\delta}
\end{align*}
A direct computation shows
\begin{align*}
	|\nabla \left( |\Psi|^{\alpha+1-n} \Psi^n \right)| \leq C \left\langle n \right\rangle |\Psi|^{\alpha}|\nabla \Psi|.
\end{align*}
Applying Proposition \ref{pro:ga1}, we have
\begin{align*}
	\left\lVert |\Psi|^{\alpha+1-n} \Psi^n \right\rVert_{\dot{H}^{1}} 
	\leq{}& C \left\langle n \right\rangle  \left\lVert \phi \right\rVert_{\infty}^{\alpha} \norm{\phi}_{\dot{H}^1} \( 1 + |\eta| \norm{\phi}_{L^{\infty}}^{\alpha} \),
\end{align*}
form which the desired follows. 

Let us move on to the case $d=2$, $3$.
Combining Lemma \ref{lem:g3a} with Proposition \ref{pro:ga1},
we see from Lemma \ref{lem_lei} and \eqref{gn:1} that
\begin{align*}
	\left\lVert  |\Psi|^{\alpha+1-n} \Psi^n \right\rVert_{\dot{H}^{\delta}} 
	\leq{}& \left\lVert \nabla \left( |\Psi|^{\alpha+1-n} \Psi^n  \right) \right\rVert_{\dot{H}^{\delta-1}} \\
	\leq{}& C \left\langle n \right\rangle^{\gamma} \left(  \norm{\phi}^{\alpha + 1-\delta}_{\infty} \left( \norm{\phi}_{\infty}^{1-\frac{1}{\delta}} \norm{\Psi}_{\dot{H}^{\delta}}^{\frac{1}{\delta}} \right)•^{\delta} + \left\lVert \phi \right\rVert_{\infty}^{\alpha} \left\lVert |\nabla|^{\delta} \Psi \right\rVert_{2} \right) \\
	\leq{}& C \left\langle n \right\rangle^{\gamma} \left\lVert \phi \right\rVert_{\infty}^{\alpha} \norm{\phi}_{H^{\delta}}\( 1 + |\eta| \norm{\phi}_{\infty}^{\alpha} \)^{\lceil \delta \rceil},
\end{align*}
as desired. This completes the proof.
\end{proof}

\section*{Acknowledgments}
M.K. was supported by JSPS KAKENHI Grant Numbers 20K14328.
H.M. was supported by JSPS KAKENHI Grant Numbers 19K14580 and 22K13941.


\end{document}